\documentclass{specgeom}
\usepackage[centertags]{amsmath}
\usepackage{amssymb,amsfonts,amsthm}
\usepackage{pictex}
\usepackage{epsfig}  %needed only if figure 2 is included as eps and not created on the fly

\theoremstyle{plain}
\newtheorem{theorem}{Theorem}[section]
\newtheorem{prop}[theorem]{Proposition}
\newtheorem{lemma}[theorem]{Lemma}
\newtheorem{cor}[theorem]{Corollary}

\theoremstyle{definition}
\newtheorem{dfn}[theorem]{Definition}

\theoremstyle{remark} 
\newtheorem{remark}[theorem]{Remark}

\theoremstyle{plain}

\numberwithin{equation}{section}

\providecommand{\mathscr}{\mathcal} % a priori mathscr is mathcal
\IfFileExists{mathrsfs.sty}{
\usepackage{mathrsfs}}         %  \mathscr produces rsfs fonts
   {%  else try euler fonts
     \IfFileExists{eucal.sty}{
        \usepackage[mathscr]{eucal}} % \mathscr produces Euscript fonts,
   {% else both unavailable, do nothing, mathscr remains mathcal
   }
}
\let\templabelenumi=\labelenumi
\renewcommand{\labelenumi}{\textup{\templabelenumi}}
\let\templabelenumii=\labelenumii
\renewcommand{\labelenumii}{\textup{\templabelenumii}}

\newcommand{\enumref}[1]{\textup{(\ref{#1})}}
\newcommand{\enum}[1]{\textup{(#1)}}

\newcommand{\plref}[1]{{\textup{\ref{#1}}}}

\newcommand{\bigsetdef}[2]{\bigl\{ #1 \,\bigm|\, #2\bigr\}}

\newcommand{\sa}{{\operatorname{sa}}}

\newcommand\SF{\operatorname{SF}}
\newcommand{\specess}{\operatorname{spec}_{\operatorname{ess}}}

\newcommand\ev{\operatorname{ev}}

\newcommand\id{\operatorname{id}}
\newcommand\im{\operatorname{im}}
\renewcommand\Im{\operatorname{Im}}
\newcommand\ind{\operatorname{ind}}

\newcommand\rank{\operatorname{rank}}

\renewcommand\Re{\operatorname{Re}}
\newcommand\sgn{\operatorname{sgn}}

\newcommand\spec{\operatorname{spec}}

\newcommand\eps{\varepsilon}

\newcommand\ga{\alpha} 
  
\newcommand\gl{\lambda}

\newcommand{\C}{\mathbb{C}}
\newcommand{\R}{\mathbb{R}}

\newcommand{\Z}{\mathbb{Z}}

\newcommand\cB{\mathscr{B}}
\newcommand\cC{\mathscr{C}}
\newcommand\cD{\mathscr{D}}

\newcommand\cF{\mathscr{F}}

\newcommand\cI{\mathscr{I}}

\newcommand\cK{\mathscr{K}}

\newcommand\cM{\mathscr{M}}
\newcommand\cN{\mathscr{N}}

\newcommand\cU{\mathscr{U}}

\newcommand{\scalar}[2]{\langle #1,#2\rangle}

\newcommand{\comment}[1]{\relax}

\newcommand{\csa}{\cC^\sa}
\newcommand{\cfsa}{\cC\cF^\sa}
\newcommand{\bsa}{\cB^\sa}
\newcommand{\bfsa}{\cB\cF^\sa}
\newcommand{\pitilde}{\widetilde\pi}

\begin{document}

\title[Uniqueness of spectral flow]{The Uniqueness of the Spectral Flow on Spaces 
of Unbounded Self--adjoint Fredholm Operators}
\author{Matthias Lesch}
\address{Mathematisches Institut,
Universit\"at zu K\"oln,
Weyertal 86--90,
50931 K\"oln,
Germany} 
\email{lesch@mi.uni-koeln.de}
\urladdr{http://www.mi.uni-koeln.de/$\sim$lesch}

% Version managemant, to be removed later
%\thanks{\fbox{Version: 1.0.3 Final revision 2004-04-19 date-processed: \today}}

%    General info
\date{January 28, 2004 and, in revised form, April 1, 2004.}
\thanks{The author was supported by Deutsche Forschungs\-gemeinschaft through the 
Sonder\-forschungs\-bereich/Transregio 12.}

\begin{abstract}
We discuss several natural metrics on spaces of unbounded
self--adjoint operators and their relations, among them the Riesz
and the graph metric. We show that
the topologies of the spaces of Fredholm operators resp. invertible
operators depend heavily on the metric. Nevertheless we prove
that in all cases the spectral flow is up to a normalization
the only integer invariant of non--closed paths which is path additive
and stable under homotopies with endpoints varying in the space of invertible
self--adjoint operators.

Furthermore we show that for certain Riesz continuous paths
of self--adjoint Fredholm operators the spectral flow can
be expressed in terms of the index of the pair of positive spectral
projections at the endpoints.

Finally we review the Cordes--Labrousse theorem on the stability of
the Fredholm index with respect to the graph metric in a modern
language and we generalize it to the Clifford index and to
the equivariant index.
\end{abstract}

\maketitle

% Remove in final version
%\tableofcontents
%\listoffigures

%\section*{Pflichtenheft}
%
%\begin{enumerate}
%\item Quote Robbin-Salamon and deal with their situation
%      \cite{RobSal:SFM}
%\item Introduction
%\item Multiplikation von bounded mit unbounded operator stetig???
%\end{enumerate}

\section{Introduction}

Let $H$ be a separable complex Hilbert space. Then it is well--known that the space 
$\bfsa=\bfsa(H)$ of bounded self--adjoint Fredholm operators has three connected
components, i.e. $\bfsa$ is the disjoint union
\begin{equation}
  \bfsa=\bfsa_+\cup \bfsa_-\cup \bfsa_*,
  \label{G-I.1}
\end{equation}
where $\bfsa_\pm$ denote the subspaces of essentially positive/negative
operators and $\bfsa_*=\bfsa\setminus\bigl(\bfsa_+\cup
\bfsa_-\bigr)$.

$\bfsa_\pm$ are trivially contractible. Atiyah and Singer \cite{AtiSin:ITS}
showed that the interesting component $\bfsa_*$ is
a classifying space for the $K^1$--functor. In particular, one has
\begin{equation}
    \begin{split}
       \pi_k(\bfsa_*,I)\simeq\begin{cases}
           \Z,& k \text{ odd},\\
           0, & k \text{ even}.
                                \end{cases}
    \end{split}
\label{G-I.2}
\end{equation}

The isomorphism
\begin{equation}
   \SF:\pi_1(\bfsa_*,I)\longrightarrow \Z
   \label{G-I.3}
\end{equation}
is the celebrated \emph{spectral flow}, although it was not addressed
as such in loc. cit. The spectral flow was introduced and generalized to
non--closed paths in the famous series of papers on spectral asymmetry
by Atiyah, Patodi, and Singer \cite{AtiPatSin:SAR}.
In the finite--dimensional context the spectral flow probably
dates even back to Morse and his index theorem.

It is impossible to give a complete account on the literature about
the spectral flow. I would like to emphasize, however, that a rigorous
definition of the spectral flow for (non--closed) continuous paths
of bounded self--adjoint Fredholm operators is non--trivial. After
the intuitively appealing approach of \cite{AtiPatSin:SAR} the spectral
flow was folklore and people did not feel the need or found it too
trivial to bother about the definition and its basic properties.
J. Phillips \cite{Phi:SAF} presented a completely different rigorous
approach to the spectral flow of bounded self--adjoint Fredholm operators.
Let us briefly summarize his definition:

\begin{dfn}\label{def-SF}
Let $f:[0,1]\to \bfsa$ be a continuous path of bounded
self--adjoint Fredholm operators. Choose a subdivision
 $0=t_0<t_1<\ldots < t_n=1$ of the interval such that there
exist $\eps_j>0, j=1,...,n$ with $\pm \eps_j\not\in\spec f(t)$ and
$[-\eps_j,\eps_j]\cap\specess f(t)=\emptyset$ for $t_{j-1}\le t\le t_j$.
Then the \emph{spectral flow} of $f$ is defined by
\[
   \SF(f):=
   \sum_{j=1}^n
     \Bigl(\rank\bigl(1_{[0,\eps_j)}(f(t_j))\bigr)-
                  \rank\bigl(1_{[0,\eps_j)}(f(t_{j-1}))\bigr)\Bigr).
\]
\end{dfn}
Here we have used the following notation which will be in effect
throughout the paper: by $1_X$ we denote the characteristic function
of $X$ and for a Borel subset $X\subset\C$ and a normal operator $T$
we denote by $1_X(T)$ the normal operator obtained by plugging
$T$ into $1_X$ via the Borel functional calculus.

It is shown in \cite{Phi:SAF} that a subdivision with the desired properties
indeed exists and that $\SF$ is well--defined, path additive, and 
homotopy invariant. Also the isomorphism \eqref{G-I.3} is reproved.

It is an easy consequence of the isomorphism \eqref{G-I.3} that
the spectral flow is up to normalization the only path additive
and homotopy invariant integer--valued function from paths of self--adjoint
Fredholm operators (see Theorem \plref{uniqueness-bounded} below
for a precise formulation).

In various branches of mathematics the spectral flow of families of
\emph{unbounded} operators arises naturally (e.g. in Floer homology,
Nicolaescu \cite{Nic:MIS}, Robbin and Salamon \cite{RobSal:SFM} 
to mention only a few). In the case
of boundary value problems one even has to deal with operators with
varying domains.
Superficially, one might be tempted to believe that the aforementioned
results for bounded operators just carry over with only minor modifications.

(Un)fortunately, this is not the case. So, denote by $\csa$ the
set of possibly unbounded self--adjoint 
operators in $H$ and by $\cfsa\subset\csa$ 
the subspace of (unbounded) self--adjoint Fredholm operators. 
At first, there exist several natural metrics on (subspaces of)
$\csa$ and results may depend on the metric. The weakest metric
is the graph or gap metric, $d_G$, which was studied systematically by Cordes
and Labrousse \cite{CorLab:IIM}. Another metric is the Riesz metric,
$d_R$, which was discussed by Nicolaescu in the unpublished note
\cite{Nic:SFO}. 

If one considers only operators with a fixed domain there is even another
metric: let $D$ be a fixed self--adjoint operator with domain
$W:=\cD(D)$. On the space $\bsa(W,H):=\bigsetdef{T\in\csa}{\cD(T)=W}$
there is another natural metric $d_W$ (see
Def. \plref{def-dW}). For $d_W$--continuous paths $f:[0,1]\to
\bsa(W,H)$ of Fredholm operators 
the spectral flow was defined by Boo{\ss}--Bavnbek and Furutani \cite{BooFur:MIF}
(with the additional assumption $f(t)-f(0)$ bounded) and
in \cite{RobSal:SFM} (with the assumption that $D$ has compact
resolvent).
Moreover in \cite[Sec. 4]{RobSal:SFM} it was shown that in their
case the spectral flow is also unique in the sense described above.

For the Riesz metric the results mentioned at the beginning of this section indeed 
carry over verbatim. Namely, in subsection \plref{secfivefour} we will
show that the natural inclusion of the pair $(\bfsa,G\bsa)$ into
$\bigl(\cfsa,G\csa,d_R\bigr)$ is indeed a homotopy
equivalence ($GX$ denotes the invertible elements in $X$). 
Hence the unbounded analogue of $\bfsa_*$, $(\cfsa_*,d_R)$ is a classifying space for the
$K^1$--functor and the analogue of \eqref{G-I.3} holds. There is a
drawback: the Riesz topology is so strong that it is hard to prove
continuity of maps into $(\csa,d_R)$. As an example consider a compact
manifold with boundary, $M$, and a Dirac operator $D$ on $M$. Elliptic
boundary conditions for $D$ are parametrized by certain
pseudodifferential projections (cf. e.g. Br\"uning and Lesch \cite{BruLes:BVPI} and
references therein) $P$ on the boundary. Denote by $D_P$ the
self--adjoint realization of $D$ with boundary condition $P$. Then
$P\mapsto D_P$ is graph continuous. Note that $P\mapsto D_P$ is a
family of operators with varying domains! See
Boo{\ss}--Bavnbek, Lesch, and Phillips \cite[Sec. 3]{BooLesPhi:UFO} for details.
It is not known, at least not to the author, whether $P\mapsto D_P$ is
Riesz continuous or not. 

In the development of the spectral flow
for paths of unbounded operators one first tried to use the Riesz
metric. But continuity proofs for simple maps like 
$\bsa\to\csa, A\mapsto D+A$ ($D$ a fixed self--adjoint operator)
are rather complicated (cf. \cite[Thm. A.8]{Phi:SFT}, {\cite[Sec. 4]{BooFur:MIF} and
Proposition \plref{S1.1}).

A drawback (if it is one!) of the weaker graph topology is that the
homotopy type $\cfsa$ is presumably more complicated. 
For graph continuous paths in $\cfsa$ it was shown in \cite{BooLesPhi:UFO}
that the definition of the spectral flow in \cite{Phi:SAF} carries
over. More importantly,
an alternative definition of the spectral flow in terms of the Cayley
transform and the classical winding number is given. This uses
results of Kirk and Lesch \cite[Sec. 6]{KirLes:EIM} which show that the spectral
flow, the Maslov index, and the classical winding number are intimately
connected.

It should come as a surprise that, as opposed to \eqref{G-I.1}, 
$(\cfsa,d_G)$ is path connected \cite[Thm. 1.10]{BooLesPhi:UFO}. In light of
this one should also be ready for surprises concerning the fundamental group.
Still, except that the spectral flow is a surjective homomorphism from
$\pi_1(\cfsa,d_G)$
onto the integers nothing about $\pi_1(\cfsa,d_G)$ is known.
Therefore, we single out the following open problem:

\newtheorem*{problem}{Problem}

\begin{problem} Find $\pi_1(\cfsa,d_G)$. Even more, is $(\cfsa,d_G)$ a
  classifying
space for the $K^1$--functor?
\end{problem}

We do not have a good guess for the answer to this problem and therefore we do
not further speculate.

Since the fundamental group of $(\cfsa,d_G)$ is not known
the uniqueness of the spectral flow on $(\cfsa,d_G)$ cannot (yet) be
proved along the lines of the case of bounded operators. 
The current paper wants, among other things, to fill this
gap and prove that as in the bounded case the 
spectral flow is up to normalization the only path additive
and homotopy invariant integer--valued function from paths in $(\cfsa,d_G)$
(see Theorem \plref{S-3.4} below for a precise formulation).

The second goal of this paper is to give an account on the various metrics
on (subspaces of) $\csa$ and their relations. It should be noted
at this point that the restriction to self--adjoint operators is not a loss
of generality as far as the space of \emph{all} unbounded operators
is concerned. 
Namely, denote by $\cC$ the set of closed densely
defined operators in $H$. Then the map
\begin{equation}
      T\mapsto \begin{pmatrix} 0 & T\\ T^* & 0\end{pmatrix}
\end{equation}
is a natural embedding $\cC(H)\hookrightarrow \csa(H\oplus H)$. So each metric
on $\csa(H\oplus H)$ 
naturally induces a metric on $\cC(H)$ and, obviously, each
metric on $\cC(H)$ induces one on $\csa(H)$. That is the reason why we restrict
ourselves to the discussion of self--adjoint operators. This approach,
admittedly,
does not cover all cases one could think of in this context: for example the
space of all closed operators with a fixed domain (for instance the Sobolev
space $H^1\subset L^2$ on a manifold) does not seem to be treatable
by this approach; see however Definition \plref{def-dW} below.

The paper is organized as follows:

In Section \plref{sectwo} we introduce the various metrics on spaces
of unbounded operators, namely the graph metric $d_G$, the Riesz
metric $d_R$, the $d_W$--metric and the norm metric $d_N$. 
The latter two are defined only
on certain subspaces of $\csa$. We show that 
$d_N\succneqq d_W\succneqq d_R\succneqq d_G$ (see Propositions \plref{S1.1} and
\plref{S1.2}), i.e.  $d_N$ is strictly stronger than $d_W$ and so on.
We put some effort in the construction of 
counterexamples which show that the metrics induce different topologies
even on relatively "small" subsets of $\csa$.

Section \plref{secthree} discusses the relation between the spectral flow
and the index of a pair of projections. More precisely, we show that
for a \emph{Riesz} continuous path $T_t$ of self--adjoint Fredholm operators
with the additional property that the domains are fixed and that
$T_t-T_0$ is compact the spectral flow is the index of the pair of positive
spectral projections at the endpoints (Theorem \plref{S-SFpair}). 
This generalizes work of Bunke \cite{Bun:SFF}
who has considered the special case of families of the form $D_t:=D+tR$
where $D$ is a Dirac operator on a compact manifold and $R$ is a
self--adjoint bundle endomorphism satisfying additional assumptions.

As an application we prove an abstract Toeplitz index theorem
(Proposition \plref{Toeplitz}).

The positive spectral projections and their relative index were used
to define the spectral flow even in the von Neumann algebra context by
Phillips \cite{Phi:SFT}.

Section \plref{secfour} presents the results of the celebrated paper by
Cordes and Labrousse \cite{CorLab:IIM} in modern language and in a very
concise form. We go slightly beyond loc. cit. and prove the stability
of the Clifford index and the $G$--index with respect to the graph metric
(Theorem \plref{stability}). Moreover we prove 
that with respect to both, $d_G$ and $d_R$, the bounded
self--adjoint operators are open and dense in $\csa$ and that $d_G$
and $d_R$ induce the norm topology on bounded operators
(Proposition  \plref{S-4.1}).

Section \plref{secfive} is the heart of the paper.  We show that 
the spectral flow can be characterized axiomatically, 
i.e. it is the only integer invariant of continuous
paths of self--adjoint Fredholm operators which satisfies
\emph{Homotopy, Concatenation,} and \emph{Normalization}. We prove
this uniqueness in 
the finite--dimensional case, the bounded case, for the graph and Riesz
metric, and for the $d_W$--metric (Theorems
\plref{uniqueness-bounded}, \plref{uniqueness-fd},
\plref{uniqueness-graph}, \plref{uniqueness-Riesz}, and
\plref{uniqueness-dW}).
On the way we generalize the method of
\cite{BooLesPhi:UFO} to show that the space of invertible elements of
$(\csa, d_G)$ is still path connected (Proposition \plref{S3.1id}).

In the Appendix we finally collect a few useful operator estimates which we
need and which did not quite fit into the course of the paper.

Part of this work was presented at the 
Workshop "Spectral geometry of manifolds with boundary and decomposition of
manifolds" held in Holbaek, Denmark, August 6--9, 2003.
I would like to thank B. Boo{\ss}--Bavnbek, G. Grubb, and K. Woj\-ciechowski for 
providing this great and stimulating environment and for compiling the
current proceedings.

Furthermore, I would like to thank B. Boo{\ss}--Bavnbek for his constant interest
in this work and for numerous discussions.

Finally, I wish to thank the anonymous referee for helpful comments.

%\newpage
\section{Topologies on spaces of unbounded self--adjoint operators}
\label{sectwo}

In this section we will discuss various natural metrics and their topologies
on spaces of unbounded operators. These are used frequently in the literature.
However, a systematic comparison does not seem to be available except
for the Riesz and the graph topology, see e.g. \cite{Nic:SFO},
\cite{BooLesPhi:UFO}. 

\subsection*{Notations}
 Let $H$ be a separable complex Hilbert space. First
let us introduce some notation for various spaces of operators in $H$:
\begin{align*}
 \cC(H) & \text{ closed densely defined operators in $H$},\\
 \cB(H) & \text{ bounded linear operators $H\to H$},\\
 \cU(H) & \text{ unitary operators $H\to H$},\\
 \cK(H) & \text{ compact linear operators $H\to H$},\\
 \cB\cF(H) & \text{ bounded Fredholm operators $H\to H$},\\
 \cC\!\cF(H) & \text{ (closed) densely defined Fredholm operators in $H$}.
\end{align*}

If no confusion is possible we will omit ``$(H)$'' and write $\cC, \cB,
\cK$ etc. By $\cC^{\sa}, \cB^{\sa}$ etc. we denote the set of self-adjoint
elements in $\cC, \cB$ etc.

\medskip
$\csa$ carries two natural metrics, the \emph{Riesz metric} and the \emph{gap}
or \emph{graph metric}. The Riesz metric is given by
\begin{equation}
    d_R(T_1,T_2):=\|F(T_1)-F(T_2)\|,
\end{equation}
where
\begin{equation}
    F(T):=T(I+T^2)^{-1/2}.
\end{equation}
The graph metric is given by 
\begin{equation}
    d_G(T_1,T_2):=\frac 12 \|\kappa(T_1)-\kappa(T_2)\|=\|(T_1+i)^{-1}-
    (T_2+i)^{-1}\|,
\end{equation}
where 
\begin{equation}
    \kappa(T)=(T-i)(T+i)^{-1}
\end{equation}
is the Cayley transform (cf. \cite[Thm. 1.1]{BooLesPhi:UFO}). For alternative
descriptions of the graph metric see Nicolaescu \cite[Appendix A]{Nic:GSG} or Kato
\cite{Kat:PTL}.

If we restrict ourselves to operators with a fixed domain then there
are even more metrics: let $D$ be a fixed self--adjoint operator in $H$. 
The domain of $D$, $W:=\cD(D)$,
equipped with the graph scalar product,
\begin{equation}
\scalar{x}{y}_D:=\scalar{x}{y}+\scalar{Dx}{Dy},
\end{equation}
is then a Hilbert space which is continuously embedded in $H$.

\begin{dfn}\label{def-dW} $\bsa(W,H):=\bigsetdef{T\in\csa}{\cD(T)=W}.$
$\bsa(W,H)$ is equipped with the metric
\[
   d_W(T_1,T_2):=\|T_1-T_2\|_{W\to H}=\|(T_1-T_2)(I+D^2)^{-1/2}\|_{H\to H}.
\]
\end{dfn}

Note that if $T_1\in\bsa(W,H)$ is invertible as an element of
$\csa$ then $T_1^{-1}$ maps $H$ continuously into $W$ and thus
for any $T_2\in\bsa(W,H)$ the operator $T_2T_1^{-1}$ is a bounded
operator $H\to H$. This will be used in the sequel without further notice.

On the subspace
\begin{equation}
    D+\bsa=\bigsetdef{D+C}{C\in\bsa}\simeq \bsa
\end{equation}
we have additionally the norm distance
\begin{equation}
    d_N(D+C_1,D+C_2)=\|C_1-C_2\|.
\end{equation}

\begin{prop}\label{S1.1}The natural maps
\[\begin{split}
   (\bsa,d_N)&\stackrel{\alpha}{\longrightarrow}
   (\bsa(W,H),d_W)\stackrel{\beta}{\longrightarrow}
   (\csa,d_R)\stackrel{\id}{\longrightarrow}
   (\csa,d_G)\\
     C&\longmapsto D+C
  \end{split}
\]
are continuous. Here $\beta$ is the natural inclusion.
\end{prop}
\begin{remark}\indent\par

\enum{1} The continuity of the
identity map $(\csa,d_R)\to (\csa,d_G)$ was observed by
Nicolaescu \cite[Lemma 1.2]{Nic:SFO}.

\medskip
\enum{2}
The continuity of $\beta$ generalizes \cite[Thm. 4.8 and
Cor. 4.9]{BooFur:MIF} where it is proved that the composition
map $\beta\circ\alpha$ is continuous.
\end{remark}

\begin{proof}
\enum{1} For $C_1,C_2\in\bsa$ we have 
\begin{equation}
    \begin{split}
     d_W(D+C_1,D+C_2)&=\|(C_1-C_2)(I+D^2)^{-1/2}\|\\
                &\le \|C_1-C_2\|\, \|(I+D^2)^{-1/2}\|\\
                &\le d_N(D+C_1,D+C_2),
    \end{split}
\end{equation}
i.e. $d_W\le d_N$ and hence $\alpha$ is continuous.

\medskip
\enum{2} For completeness we briefly recall Nicolaescu's \cite{Nic:SFO} argument
to prove the continuity of the identity map 
$(\csa,d_R)\stackrel{\id}{\longrightarrow} (\csa,d_G)$: for
$T\in\csa$ we have
\begin{equation}
   (T+i)^{-1}=(T-i)(I+T^2)^{-1}=(I+T^2)^{-1/2}F(T)-i(I+T^2)^{-1}
\end{equation}
and
\begin{equation}\label{G1.9}
   (I+T^2)^{-1}=I-F(T)^2.
\end{equation}
Hence, if $F(T_n)\to F(T)$ then $(I+T_n^2)^{-1}\to (I+T^2)^{-1}$ and
thus also $(I+T_n^2)^{-1/2}\to (I+T^2)^{-1/2}$.
Consequently $(T_n+i)^{-1}\to (T+i)^{-1}$.

\medskip
\enum{3} The continuity of $\beta$ is more complicated.
We fix a $T\in\bsa(W,H)$ and we have to prove the continuity of 
$F$ at $T$. Put
\begin{equation}
     M:=\|(T\pm i)^{-1}(I+D^2)^{1/2}\|=\|(I+D^2)^{1/2}(T\pm i)^{-1}\|.
\end{equation}
Let $0<q<\frac 12$ and consider $\widetilde T\in\bsa(W,H)$ with
$d_W(T,\widetilde T)\le \frac qM$. Then we have
\begin{equation}\label{G1.11}
     \|(T-\widetilde T)(T\pm i)^{-1}\|, \|(T\pm i)^{-1}(T-\widetilde T)\|\le q.
\end{equation}
The Neumann series then immediately implies
\begin{equation}
    \|(\widetilde T+i)^{-1}(T+i)\|\le \frac{1}{1-q}.
\end{equation}
Thus, for $x\in H$ we have
\begin{equation}
     \|(\widetilde T+i)^{-1}x\|\le \frac{1}{1-q}\|(T+i)^{-1}x\|
\end{equation}
and
\begin{equation}
     \|(T+i)^{-1}x\|\le \|(T+i)^{-1}(\widetilde T+i)\|\,\|(\widetilde
     T+i)^{-1}x\| \le (1+q)\|(\widetilde T+i)^{-1}x\|.
\end{equation}
This implies the operator inequalities
\begin{equation}
    \frac{1}{(1+q)^2}|T+i|^{-2}\le |\widetilde T+i|^{-2}\le
    \frac{1}{(1-q)^2}|T+i|^{-2}.
\end{equation}
Since the square root is an operator--monotonic increasing function 
(Kadison and Ringrose \cite[Prop. 4.2.8]{KadRin:FTOI}) we may take the square root of these
inequalities and after subtracting $|T+i|^{-1}$ we arrive at
\begin{equation}
    -\frac{q}{1+q}|T+i|^{-1}\le |\widetilde T+i|^{-1}-|T+i|^{-1}\le
    \frac{q}{1-q}|T+i|^{-1}.
\end{equation}
This gives
\begin{equation}
   \|\,|T+i|^{1/2}|\widetilde T+i|^{-1}|T+i|^{1/2}-I\|\le \frac{q}{1-q}.
\end{equation}
In the following series of estimates we are going to use the estimate 
Proposition \plref{S-A.1} several times:
\begin{equation}\begin{split}
         \|F(T)&-F(\widetilde T)\|\\
         &\le \bigl\|\,|i+T|^{-1/2}(F(T)-F(\widetilde T))
     |i+T|^{1/2}\|\\
         &\le \bigl\|\,|i+T|^{-1/2}(T-\widetilde T)|i+T|^{-1/2}\bigr\|\\
            &\quad  +\bigl \|\, |T+i|^{-1/2}\bigl(\widetilde T(|i+T|^{-1}-
               |i+\widetilde T|^{-1})\bigr)|i+T|^{1/2}\bigr\|\\
        &\le \||i+T|^{-1}(T-\widetilde T)\|\\
           &\quad +\|\,|i+T|^{-1/2}\widetilde T |i+T|^{-1/2}\|\,\|I-|i+T|^{1/2}|i+
                    \widetilde T|^{-1}|i+T|^{1/2}\|\\
        &\le q+\|\,|i+T|^{-1}\widetilde T\| \frac{q}{1-q}\\
        &\le q(1+\frac{1+q}{1-q}).
                \end{split}
\end{equation}
This shows that if $d_W(T_n,T)\to 0$ then $F(T_n)\to F(T)$ and we are done.
\end{proof}

By a famous example due to Fuglede (\cite[Rem. 1.5]{Nic:SFO},
\cite[Ex. 2.14]{BooLesPhi:UFO}) the Riesz topology on $\csa$
is strictly stronger than the graph topology. 
The counterexamples in loc. cit. even have fixed domain, i.e.
a sequence of the form $T_n=D+C_n, C_n\in\bsa,$ is constructed such
that $T_n$ converges in the graph but not in the Riesz topology.
%However, in the Fuglede example 
%$(C_n)$ is not a bounded sequence, i.e. $\lim\limits_{n\to\infty}
%\|C_n\|=\infty.$
We will refine the Fuglede example and show that
the four topologies induced by $d_N, d_W, d_R, d_G$ are all different.

Before let us introduce a bit of notation. 
For metrics $d_1,d_2$ on a metric space we write $d_1\succeq d_2$
($d_1\succneqq d_2$) if the topology induced by $d_1$ 
is (strictly) stronger than the one induced by $d_2$.
Of course, if $d_1\ge d_2$ then $d_1\succeq d_2$ but the converse
need not be true.

\begin{prop}\label{S1.2}
Let $H$ be a separable complex Hilbert space and
let $D$ be a self--adjoint operator in $H$ with compact resolvent.
\begin{enumerate}
\item\label{S1.2.1} On $D+\bsa$ we have $d_N\succneqq d_W\succneqq d_R\succneqq d_G$.
\item\label{S1.2.2} For fixed $R\ge 0$ we have on the space 
\[\bigsetdef{D+C}{C\in\bsa,\, \bigl\||D+i|^{-1}C|D+i|\bigr\|\le R}\]
that $d_W\succeq d_R$ and $d_R\succeq d_W$, i.e.
$d_W$ and $d_R$ induce the same topology on this subset of $D+\bsa$.
\end{enumerate}
\end{prop}
\begin{proof} It follows from Proposition  \plref{S1.1} that 
$d_N\succeq d_W\succeq d_R\succeq d_G$.

Next we prove that on $\bigsetdef{D+C}{C\in\bsa,\,
\bigl\||D+i|^{-1}C|D+i|\bigr\|\le R}$ we also have $d_R\succeq d_W$. 

%Old version: wrong!!
\comment{For $C\in\bsa, \|C\|\le R$, we have
\begin{equation}\begin{split}
      \bigl\|(&I+(D+C)^2)^{1/2}(I+D^2)^{-1/2}\bigr\|\\
     &\le \bigl\|(I+(D+C)^2)^{1/2}(D+C+i)^{-1}\bigr\|\\
         &\quad
        \bigl\|(D+C+i)(D+i)^{-1}\bigr\|\,\bigl\|(D+i)(I+D^2)^{-1/2}\bigr\|\\
    &\le \bigl\|(D+C+i)(D+i)^{-1}\bigr\|\\
    &\le 1+R.
                \end{split}
\end{equation}
In the third inequality we have used
that $(I+(D+C)^2)^{1/2}(D+C+i)^{-1}$ is unitary.

Thus we have for $C,\widetilde C\in\bsa$, $\|C\|, \|\widetilde C\|\le R$.
\begin{equation}
    \begin{split}
                d_W(&D+C,D+\widetilde C)\\
            &= \|(C-\widetilde C)(I+D^2)^{-1/2}\|\\
               &\le (1+R)\|(C-\widetilde C)(I+(D+C)^2)^{-1/2}\|\\
        &\le (1+R)\Bigl(\|F(C)-F(\widetilde C)\|\\
         &\quad +\|\widetilde C\bigl((I+(D+C)^2)^{-1/2}-(I+(D+\widetilde
    C)^2)^{-1/2}\bigr)\|\Bigr)\\
       &\le (1+R)\Bigl(\|F(C)-F(\widetilde C)\|\\
           &\quad +R\|(I+(D+C)^2)^{-1/2}-(I+(D+\widetilde
    C)^2)^{-1/2}\|\Bigr).
    \end{split}
\end{equation}
If $d_R(D+C_n,D+C)\to 0$ then (cf. the argument after \eqref{G1.9})
$(I+(D+C_n)^2)^{-1/2}\to (I+(D+ C)^2)^{-1/2}$ 
and hence $d_W(D+C_n,D+C)\to 0$.
}%end of wrong version

Let $C_n,C\in\bsa$, $\||D+i|^{-1}C_n|D+i|\|\le R, \||D+i|^{-1}C|D+i|\|\le R$, and assume that
$d_R(D+C_n,D+C)\to 0$, i.e. $F(D+C_n)\to F(D+C)$, $n\to\infty$.

We note that it follows from Proposition \plref{S-A.1} that the operators
\begin{equation}
|D+i|C_n|D+i|^{-1},\quad |D+i|C|D+i|^{-1}
\end{equation}
are bounded (and defined on all of $H$) and satisfy the same norm bound.

Consider the identity
\begin{equation}\begin{split}
    &F(D+C_n)-F(D+C)\\
     &=(D+C_n)\bigl[|D+C_n+i|^{-1}-|D+C+i|^{-1}\bigr]+(C_n-C)|D+C+i|^{-1}.
                \end{split}\label{G1.23}
\end{equation}
We have to show that $\|(C_n-C)|D+i|^{-1}\|\to 0$. We first note that
it suffices to show that $(C_n-C)|D+C+i|^{-1}\to 0$ strongly. Indeed,
if this is the case then for $x\in\cD(D)$ we have
$(C_n-C)x=\Bigl((C_n-C)(D+C+i)^{-1}\Bigr)\Bigl((D+C+i)x\Bigr)\to 0$. Hence
$(C_n-C)\to 0$ strongly on the dense subspace $\cD(D)$. Since in view of
Proposition \ref{S-A.1} $\|C_n-C\|\le R$ is uniformly bounded we infer that
$(C_n-C)\to 0$ strongly on $H$. Now since $D$ has compact resolvent
$|D+i|^{-1}$ is compact and since multiplication from the right by compact
operators turns strongly convergent sequences into uniformly convergent
sequences we indeed conclude that $\|(C_n-C)|D+i|^{-1}\|\to 0$.

To prove that $(C_n-C)|D+C+i|^{-1}\to 0$ strongly we assume the contrary.
Then there is an $x\in H$ and an $\eps>0$ such that after possibly considering
a subsequence we have
\begin{equation}
    \Bigl\|(C_n-C)|D+C+i|^{-1}x\Bigr\|\ge \eps.\label{G1.24}
\end{equation}
Again, since $\|C_n-C\|$ is uniformly bounded and since $\cD(D)$ is dense
in $H$ we may assume that $x\in\cD(D)$.

Since $F(D+C_n)\to F(D+C)$ (cf. the argument after \eqref{G1.9}) we  have
\begin{equation}
     x_n:=\bigl[|D+C_n+i|^{-1}-|D+C+i|^{-1}\bigr]x\to 0,\quad n\to\infty.
\end{equation}
Applying Proposition \ref{S-A.2} with $\alpha=2,\beta=-1$ (and $\alpha=0, \beta=-1$
and repeatedly using the boundedness of $\||D+i|C_n|D+i|^{-1}\|, 
\||D+i|C|D+i|^{-1}\|$)
we infer that
\begin{equation}
   y_n=(D+C_n)x_n=(D+C_n)\bigl[|D+C_n+i|^{-1}-|D+C+i|^{-1}\bigr]x
\end{equation}
is a bounded sequence in $\cD(D)$. Since $D$ has compact resolvent
the inclusion $\cD(D)\hookrightarrow H$ is compact and thus a subsequence
of $(y_n)$ converges in $H$.

Summing up we have proved that there is a subsequence
$x_{n_k}$ such that $x_{n_k}\to 0$ and such that
(since $C_n$ is bounded) $Dx_{n_k}$ converges in
$H$. But $D$ is a closed operator, hence $Dx_{n_k}\to 0$ and thus 
$y_{n_k}\to 0$.

Plugging $x$ into the identity \eqref{G1.23} we arrive at 
$(C_{n_k}-C)|D+C+i|^{-1}x\to
0$ contradicting \eqref{G1.24}.

\bigskip
Finally we are going to present three counterexamples which prove
the claimed $\succneqq$ relations:

Since $D$ has compact resolvent there is an orthonormal basis
$(e_k)_{k=1}^\infty$ of eigenvectors, $De_k=\lambda_k e_k$, and
$\lim\limits_{k\to\infty}|\gl_k|=\infty.$

\medskip
\enum{1} Let $C_n\in \bsa$ be defined by
\begin{equation}
      C_n e_k:=\begin{cases} e_n, &k=n,\\ 0, & \text{otherwise}.\end{cases}
\end{equation}
Then $C_n$ is a self--adjoint rank--one operator, $\|C_n\|=1$, and hence
$d_N(D+C_n,D)=1$. On the other hand, however, we find
\begin{equation}
      C_n (I+D^2)^{-1/2}e_k=\begin{cases} (1+\gl_n^2)^{-1/2}e_n, &k=n,\\ 0, & \text{otherwise},\end{cases}
\end{equation}
and thus
\begin{equation}
     d_W(D+C_n,D)=(1+\gl_n^2)^{-1/2}\xrightarrow{n\to\infty}0.
\end{equation}

This proves $d_N\succneqq d_W$ in part (1) and (2) of the Proposition.

\medskip
\enum{2} Next we put
\begin{equation}
      C_n e_k:=\begin{cases} \gl_ne_n, &k=n,\\ 0, & \text{otherwise}.\end{cases}
\end{equation}
Again, $C_n$ is a self--adjoint rank--one operator, $\|C_n\|=|\gl_n|$, and 
\begin{equation}\begin{split}
      d_W(D+C_n,D)&=\|C_n (I+D^2)^{-1/2}\|\\
                  &\ge\|C_n(I+D^2)^{-1/2}e_n\|\\
                  &=\frac{|\gl_n|}{\sqrt{1+\gl_n^2}}\xrightarrow{n\to\infty} 1.
                \end{split}
\end{equation}
On the other hand, however, we find
\begin{equation}
      \|F(D+C_n)-F(D)\|=\Bigl|\frac{2\gl_n}{\sqrt{1+(2\gl_n)^2}}-\frac{\gl_n}{\sqrt{1+\gl_n^2}}
                         \Bigr|  \xrightarrow{n\to\infty} 0.
\end{equation}

This proves $d_W\succneqq d_R$ on $D+\bsa$.

\medskip
\enum{3} The following is the famous example due to Fuglede
(\cite[Rem. 1.5]{Nic:SFO}, \cite[Ex. 2.14]{BooLesPhi:UFO}): put
\begin{equation}
      C_n e_k:=\begin{cases} -2\gl_ne_n, &k=n,\\ 0, & \text{otherwise}.\end{cases}
\end{equation}
Then
\begin{equation}
  d_G(D+C_n,D)=\Bigl|(-\gl_n+i)^{-1}-(\gl_n+i)^{-1}\Bigr|\to 0,\quad n\to\infty.
\end{equation}
On the other hand, however,
\begin{equation}\begin{split}
  \|F(D+&C_n)-F(D)\|\\
&\ge \bigl\|\bigl(F(D+C_n)-F(D)\bigr)e_n\bigr\|=|2\gl_n(1+\gl_n^2)^{-1/2}|\to
  2,\quad n\to\infty,
                \end{split}
\end{equation}
and hence $d_R\succneqq d_G$ on $D+\bsa$.
\end{proof}
\begin{remark}\indent\par

\enum{1} By a result due to Nicolaescu \cite[Prop. 1.4]{Nic:SFO} Riesz convergence
can also be characterized as follows: let $f:\R\to\C$ be any continuous
function with $f(x)=1$ for $x>>1$ and $f(x)=-1$ for $x<<-1$. Then a sequence
$T_n\in\csa$ is $d_R$--convergent if and only if it is $d_G$--convergent
\emph{and} $f(T_n)$ is convergent.

In particular this implies that if $(T_n)$ is a sequence of operators
with $T_n\ge -C$ for some fixed $C$ then $(T_n)$
is $d_R$--convergent if and only if it is $d_G$--convergent.

\medskip
\enum{2} Propostion \plref{S1.2} \enumref{S1.2.2}
is sharp in the sense that in general $d_W\succneqq d_R$
even on $\bigsetdef{D+C}{C\in\bsa,\,\bigl\|C\bigr\|\le R}$.
To see this consider
\begin{equation}
      C_n e_k:=\begin{cases} e_n, &k=1,\\ e_1, & k=n,\\ 0,&\text{otherwise}.\end{cases}
\end{equation}
$C_n$ is a self--adjoint rank--two operator, $\|C_n\|=1$.
We have
\begin{equation}\|(D+C_n+i)^{-1}(D+i)\|\le 1+\|C_n\|=2\end{equation}
and thus
\begin{equation}\begin{split}
    d_G(D+C_n,D)&=\|(D+C_n+i)^{-1}-(D+i)^{-1}\|\\
                &=\|(D+C_n+i)^{-1}C_n(D+i)^{-1}\|\\
                &\le 2\|(D+i)^{-1}C_n(D+i)^{-1}\|.
                \end{split}
\end{equation}
Furthermore,
\begin{equation}
   (D+i)^{-1}C_n(D+i)^{-1}e_k=\begin{cases} \frac{1}{(\gl_1+i)(\gl_n+i)}e_n,&  k=1,\\
                                            \frac{1}{(\gl_1+i)(\gl_n+i)}e_1, &
                                            k=n,\\  
                                             0, & \text{otherwise}.
                              \end{cases}
\end{equation}
Consequently
\begin{equation}
    d_G(D+C_n,D)\le \frac{2}{\sqrt{1+\gl_n^2}}\xrightarrow{n\to\infty} 0.
\end{equation}
With a little more effort one can show that also $d_R(D+C_n,D)\to 0$. However,
if e.g. $D$ is essentially positive then $d_R(D+C_n,D)\to 0$ follows already 
from the previous remark.

On the other hand, however,
\begin{equation}
       d_W(D+C_n,D)=\|C_n(I+D^2)^{-1/2}\|\ge
       \|C_n(I+D^2)^{-1/2}e_1\|=\frac{1}{\sqrt{1+\gl_1^2}},
\end{equation}
and thus $D+C_n$ does not converge to $D$ in the $d_W$--metric. 

In view of Proposition  \plref{S1.2} \enumref{S1.2.2} this means that
$\|(D+i)^{-1}C_n(D+i)\|$ must be unbounded. Indeed, 
\begin{equation}\begin{split}
   \|(D+i)^{-1}C_n(D+i)\|&\ge\|(D+i)^{-1}C_n(D+i)e_n\|\\
                         &=|\gl_n+i||\gl_1+i|^{-1}\to \infty,\quad n\to\infty.
\end{split}
\end{equation}

\medskip
\enum{3} We leave it as an intriguing open problem to find out whether
the metrics $d_R$ and $d_G$ induce equivalent topologies on 
$\bigsetdef{D+C}{C\in\bsa,\,\bigl\|C\bigr\|\le R}$.
\end{remark}

\section{The spectral flow, index of a pair of projections, and
the abstract Toeplitz index theorem}
\label{secthree}

In this section we relate the spectral flow of certain Riesz
continuous paths to the index of the positive spectral projections
at the endpoints. This generalizes the work of Bunke \cite{Bun:SFF}
who has considered the special case of families of the form $D_t:=D+tR$
where $D$ is a Dirac operator on a compact manifold and $R$ is a
self--adjoint bundle endomorphism satisfying additional assumptions.

\begin{dfn}\label{S4.1} Let $P,Q$ be orthogonal projections
in the Hilbert space $H$. The pair $(P,Q)$ is called a Fredholm pair
if the map $Q:\im P\to\im Q$ is a Fredholm operator. The index
of this operator is denoted by $\ind(P,Q)$.
\end{dfn}

As pointed out by the referee the notion of the index of a pair of projections
was introduced by Brown, Douglas, and Fillmore \cite{BroDouFil:UEM} who
called it the "essential codimension". Boo{\ss}--Bavnbek and Wojciechowski
\cite[p. 129 ff]{BooWoj:EBP} used the terminology "virtual codimension".

Avron, Seiler, and Simon
\cite{AvrSeiSim:IPP} gave a systematic account of Fredholm pairs. In particular
they showed that a pair $(P,Q)$ of orthogonal projections is Fredholm
if and only if $\pm 1\not\in\specess (P-Q)$ \cite[Prop. 3.1]{AvrSeiSim:IPP}. 
The latter means that the images
$\pi(P), \pi(Q)$ in the Calkin algebra $\cB/\cK$ satisfy 
\begin{equation}
\|\pi(P)-\pi(Q)\|<1.
\end{equation}
In \cite{Phi:SFT} it was shown that the index of a pair of projections 
can be developed solely
from this inequality and that it generalizes to arbitrary semifinite von
Neumann factors. Furthermore, this was used to give a completely 
general definition of spectral flow
for continuous paths in the bounded case and hence in the Riesz metric as
well.

For the basic properties of Fredholm pairs we refer to \cite{AvrSeiSim:IPP}
and \cite{BruLes:BVPI} whose results we use freely. 
We only record the following which
is proved in \cite[Lemma 2.4]{Bun:SFF} only in a special case.

\begin{lemma}\label{S4.2} Let $(P(t), Q(t)), 0\le t\le 1$, 
be a norm continuous
path of Fredholm pairs. Then $\ind(P(0),Q(0))=\ind(P(1),Q(1))$.
\end{lemma}
\begin{proof} The proof follows the one in
\cite[Lemma 2.4]{Bun:SFF}. As in loc. cit. we emphasize that
the result is not standard since domain and range of $Q(t):\im P(t)\to
\im Q(t)$ varies with $t$.

By a standard fact often used in operator K-theory 
(Blackadar \cite[Prop. 4.3.3]{Bla:KTO})
there exist continuous families of unitaries
$U,V:[0,1]\to \cU, U(0)=V(0)=I$ such that $P(t)=U(t)P(0)U(t)^*$ and
$Q(t)=V(t)Q(0)V(t)^*$.
Hence
\begin{equation}
  \ind(P(t),Q(t))=\ind\Bigl(Q(0)V(t)^*U(t)P(0):\im(P(0))\to \im (Q(0))\Bigr).
\end{equation}
Now $Q(0)V(t)^*U(t)P(0)$ is a norm--continuous family of Fredholm
operators between fixed Hilbert spaces. Thus the index does
not depend on $t$ as claimed.
\end{proof}

\begin{lemma}\label{S4.3} Let $f:[0,1]\to (\cfsa,d_R)$ be a Riesz continuous
path of self--adjoint Fredholm operators. Furthermore, assume
that $\gl\not\in\spec f(t)$ for all $t$. Then the path of spectral
projections $t\mapsto 1_{(\gl,\infty)}(f(t))$ is norm--continuous.
\end{lemma}

In view of the Fuglede example (see \enum{3} in the proof of Proposition \plref{S1.2})
we cannot expect this to hold for graph continuous paths.

\begin{proof} We first note that the Riesz transform of $f$, $F(f(t))$,
is a norm--continuous path of bounded self--adjoint Fredholm operators and
by the Spectral Theorem we have
\begin{equation}
  1_{(\gl,\infty)}(f(t))=1_{(F(\gl),\infty)}\bigl(F(f(t))\bigr).
\end{equation}
Hence we are reduced to the case of a norm--continuous family of
bounded operators for which the claim is (fairly) clear in view of
the functional calculus. Namely, since 
\begin{equation}
      \spec\bigl(F(f(t))\bigr)\subset [-1,1]
\end{equation}
we have
\begin{equation}\label{G4.3}
    1_{(F(\gl),\infty)}(F(f(t)))=\frac{1}{2\pi i}\oint_{|z-(F(\gl)+2)|=2}
          \Bigl(z-F(f(t))\Bigr)^{-1}dz.
\end{equation}
Now the right hand side of \eqref{G4.3} depends continuously on $t$.
\end{proof}

Let $T\in\csa$. We recall that a symmetric operator $S$ with $\cD(S)\supset
\cD(T)$ is called $T$--compact 
(Kato \cite[Sec. IV.1.3]{Kat:PTL}) if $S(T+i)^{-1}$ is a compact operator.
Note that in this case $T+S$ with domain $\cD(T)$ is a self--adjoint
operator, too.

\begin{prop} Let $T\in\csa$ and let $S$ be a $T$--compact symmetric
operator in $H$. Then the difference of the Riesz transforms,
\[  F(T+S)-F(T),\]
is compact.
\end{prop}
The proof of this intuitively clear result is more complicated than
expected. 

\begin{proof}\enum{1} We first deal with a special case: assume
for the moment that $S$ is bounded, compact and $\im S\subset\cD(T)$.
%Then, by the Closed Graph Theorem $S$ is in fact also a bounded operator
%from $H$ into $\cD(T)$. 
Hence $R:=(T+S)^2-T^2=TS+S(T+S)$ is defined on $\cD(T)=\cD(T+S)$.
We have 
\begin{equation}
       F(T+S)-F(T)=(T+S)\bigl(|T+S+i|^{-1}-|T+i|^{-1}\bigr)+
                            S|T+i|^{-1}.
\end{equation}
In view of the assumptions on $S$ the operators 
$S\bigl(|T+S+i|^{-1}-|T+i|^{-1}\bigr)$ and $S|T+i|^{-1}$ are compact.
It remains to prove that $T\bigl(|T+S+i|^{-1}-|T+i|^{-1}\bigr)$ is compact.
Using the resolvent equation we find (cf. Remark \plref{S-A.3} and \eqref{G-A.8})
\begin{equation}\begin{split}
    T\bigl( &|T+S+i|^{-1}-|T+i|^{-1}\bigr)\\
     &=-\frac{2}{\pi}\int_0^\infty
      T(I+T^2+x^2)^{-1}(TS+S(T+S))(I+(T+S)^2+x^2)^{-1}dx
                \end{split}\label{G4.6}
\end{equation}
Now the Spectral Theorem gives the estimates
\begin{equation}
     \|T^r(I+T^2+x^2)^{-1}T^s\|=O(x^{-2+r+s}), \quad x\to\infty,
     \qquad r,s\in\{0,1\},
\end{equation}
and similarly for $T+S$ in place of $T$.
This shows that the integrand in \eqref{G4.6} is a continuous function
with values in the compact operators which is $O(x^{-2})$ as $x\to \infty$.
Hence the integral in \eqref{G4.6} is a compact operator, too.

\enum{2} Treating the general case we introduce the spectral 
projections $P_n:=1_{[-n,n]}(T)$ and put $S_n:=P_nSP_n$. Since
$P_n$ maps $H$ continuously into $\cD(T)$ we find that $S_n$ is a
bounded compact operator with $\im(S_n)\subset\cD(T)$, hence the
situation \enum{1} applies to $S_n$. Now consider
\begin{equation}
       (S_n-S)(T+i)^{-1}=(P_n-I)S(T+i)^{-1}P_n+S(T+i)^{-1}(P_n-I).
\end{equation}
$P_n$ converges to $I$ strongly and is norm bounded by $1$. Since
$S(T+i)^{-1}$ is compact we find that $(S_n-S)(T+i)^{-1}$ converges
in the norm to $0$. In other words $T+S_n$ converges to $T+S$ in the
$d_{\cD(T)}$--metric. Then, in view of Proposition  \plref{S1.1},
$T+S_n$ converges to $T+S$ also in the Riesz metric. 
Consequently $F(T+S_n)-F(T)$
converges to $F(T+S)-F(T)$.

$F(T+S_n)-F(T)$ is compact by the proved case \enum{1} and thus
$F(T+S)-F(T)$ is compact, too.
\end{proof}

\begin{cor}\label{S-compactness} Under the assumptions of the previous Proposition
let $\gl\not\in\specess T$. Then
the difference of the spectral projections
$1_{[\gl,\infty)}(T+S)-1_{[\gl,\infty)}(T)$ is a compact operator.
\end{cor}
\begin{proof} In light of the previous Proposition it suffices to
prove the claim for bounded $T$ and compact $S$. Otherwise replace
$T$ by $F(T)$ and $S$ by $F(T+S)-F(T)$.

Since $S$ is compact we have $\specess(T)=\specess(T+S)$. Hence
$\gl$ is at most an isolated 
eigenvalue of finite multiplicity of $T$ or $T+S$.
Thus we may choose $\mu<\gl$ such that $\spec(T)\cap (\mu,\gl)=
\spec(T+S)\cap (\mu,\gl)=\emptyset$. Then
\begin{equation}
     1_{[\gl,\infty)}(T)=1_{(\mu,\infty)}(T),\quad 
     1_{[\gl,\infty)}(T+S)=1_{(\mu,\infty)}(T+S).
\end{equation}

Now choose $a>\mu$ such
that $\sup\bigl(\spec(T)\cup\spec(T+S)\bigr)<2a-\mu$. 
Then
\begin{equation}\begin{split}
     1_{[\gl,\infty)}(T&+S)-1_{[\gl,\infty)}(T)\\
     &=\frac{1}{2\pi i} \oint_{|z-a|=a-\mu} (z-T-S)^{-1}-(z-T)^{-1}dz\\
     &=\frac{1}{2\pi i} \oint_{|z-a|=a-\mu} (z-T)^{-1}S(z-T-S)^{-1}dz,
                \end{split}
\end{equation}
and this is compact since $S$ is compact.
\end{proof}

\begin{theorem}\label{S-SFpair} Let $[0,1]\ni t\mapsto T_t\in (\cfsa,d_R)$
be a Riesz continuous path of self--adjoint Fredholm operators.
Assume furthermore that the domain of $T_t$ does not depend on $t$,
$\cD(T_t)=\cD(T_0)$, and that for $t\in [0,1]$ the difference
$T_t-T_0$ is $T_0$--compact.
Then the pair $\bigl(1_{[0,\infty)}(T_1),1_{[0,\infty)}(T_0)\bigr)$ 
is Fredholm and
\[ \SF(T_t)_{t\in [0,1]}=\ind\bigl(1_{[0,\infty)}(T_1),1_{[0,\infty)}(T_0)\bigr).\]
\end{theorem}

We single out a special case which will be of interest in the proof
of the uniqueness of the spectral flow.

\begin{cor}\label{S-3.7} Let $(T_t)_{t\in [0,1]}$ be a continuous path of self--adjoint
complex $n\times n$ matrices. Then
\[
   \SF(T_t)_{t\in [0,1]}=
   \rank\bigl(1_{[0,\infty)}(T_1)\bigr)-\rank\bigl(1_{[0,\infty)}(T_0)\bigr).
\]
\end{cor}
\begin{proof} For orthogonal projections $P,Q$ in a finite--dimensional
Hilbert space we clearly have
$\ind(P,Q)= \rank P-\rank Q.$\end{proof}

\begin{proof}[Proof of Theorem \plref{S-SFpair}] 
We may choose a subdivision $0=t_0<t_1<\ldots < t_n=1$ such that there
exist $\eps_j>0$ with $\pm\eps_j\not\in\spec(T_t)$ and
\begin{equation}
    \specess(T_t)\cap [-\eps_j,\eps_j]=\emptyset,\quad t_{j-1}\le t\le t_j,
    \quad j=1,\ldots,n.
\end{equation}
Then we have by Definition \plref{def-SF}
\begin{equation}\label{G4.12}
   \SF\Bigl((T_t)_{t\in [0,1]}\Bigr)=
     \sum_{j=1}^n \Bigl(\rank\bigl(1_{[0,\eps_j)}(T_{t_j})\bigr)-
                  \rank\bigl(1_{[0,\eps_j)}(T_{t_{j-1}})\bigr)\Bigr).
\end{equation}
In view of Lemma \plref{S4.3} we may, after refining the subdivison,
assume that for $t,t'\in [t_{j-1},t_j]$ we have
\begin{equation}
     \|1_{(\eps_j,\infty)}(T_{t'})-1_{(\eps_j,\infty)}(T_{t})\|<1.
\end{equation}
Then $1_{(\eps_j,\infty)}(T_{t'})$ maps 
$\im 1_{(\eps_j,\infty)}(T_{t})$ bijectively onto 
$\im 1_{(\eps_j,\infty)}(T_{t'})$.
Hence 
\begin{equation}\begin{split}
    \ind\bigl(1_{(\eps_j,\infty)}(T_{t_{j-1}})&:\im 1_{[0,\infty)}(T_{t_j})
    \to\im 1_{[0,\infty)}(T_{t_{j-1}}) \bigr)\\
    &=\rank\bigl(
      1_{[0,\eps_j)}(T_{t_j})\bigr)-\rank\bigl(1_{[0,\eps_j)}(T_{t_{j-1}})\bigr).
        \end{split}
\end{equation}
Furthermore, since $1_{[\eps_j,\infty)}(T_{t_{j-1}})-1_{[0,\infty)}(T_{t_{j-1}})$
is of finite rank we find
\begin{equation}\begin{split}
   \ind\bigl(&1_{[0,\infty)}(T_{t_j}),1_{[0,\infty)}(T_{t_{j-1}})\bigr)\\
   &=\rank\bigl(1_{[0,\eps_j)}(T_{t_j})\bigr)-\rank\bigl(1_{[0,\eps_j)}(T_{t_{j-1}})\bigr).\end{split}\label{G4.15}
\end{equation}
Equations  \eqref{G4.12} and \eqref{G4.15} give
\begin{equation}\label{G4.16}
    \SF\Bigl((T_t)_{t\in [0,1]}\Bigr)=\sum_{j=1}^n  
       \ind\bigl(1_{[0,\infty)}(T_{t_j}),1_{[0,\infty)}(T_{t_{j-1}})\bigr).
\end{equation}

So far we have not used the assumption that the domain of $T_t$ is independent
of $t$ and the difference $T_t-T_0$ is $T_0$--compact. Hence \eqref{G4.16} holds
for any Riesz continuous path in $\cfsa$.

Now in view of our compactness assumption Corollary \plref{S-compactness} implies that
the differences 
$1_{[0,\infty)}(T_{t_j})-1_{[0,\infty)}(T_{t_{j-1}})$
are compact. In particular the difference 
$1_{[0,\infty)}(T_1)-1_{[0,\infty)}(T_0)$ is compact and hence
$\bigl(1_{[0,\infty)}(T_1),1_{[0,\infty)}(T_0)\bigr)$ is a Fredholm pair
\cite[Thm. 3.4]{AvrSeiSim:IPP}.

Since the index of a pair of projections satisfies $\ind(P,R)=\ind(P,Q)+\ind(Q,R)$
if $P-Q$ or $Q-R$ is compact \cite[Thm. 3.4]{AvrSeiSim:IPP} the right hand
side of \eqref{G4.16} indeed equals $\ind\bigl(1_{[0,\infty)}(T_1),1_{[0,\infty)}(T_0)\bigr)$
and the Theorem is proved.
\end{proof}

We record explicitly that, as noted in the proof, equation \eqref{G4.16} holds
for any Riesz continuous path. For norm continuous paths of bounded
self--adjoint Fredholm operators this was already shown in \cite{Phi:SFT}.

\begin{cor}  Let $[0,1]\ni t\mapsto T_t\in (\cfsa,d_R)$
be a Riesz continuous path of self--adjoint Fredholm operators. Choose a
subdivision as in the beginning of the proof of Theorem \plref{S-SFpair} which
is fine enough such that for $t,t'\in [t_{j-1},t_j]$ we have
\begin{equation}
     \|1_{(\eps_j,\infty)}(T_{t'})-1_{(\eps_j,\infty)}(T_{t})\|<1.
\end{equation}
Then we have
\begin{equation}\label{G4.16copy}
    \SF\Bigl((T_t)_{t\in [0,1]}\Bigr)=\sum_{j=1}^n  
       \ind\bigl(1_{[0,\infty)}(T_{t_j}),1_{[0,\infty)}(T_{t_{j-1}})\bigr).
\end{equation}
\end{cor}

As a consequence of Theorem \plref{S-SFpair} we note the abstract Toeplitz Index Theorem
(cf. \cite[Prop. 3.1]{Bun:SFF}).

\begin{prop}\label{Toeplitz}
Let $D\in\cfsa$ and let $W\in\cU$ be a unitary operator with
$W^*(\cD(D))=\cD(D)$ and $[D,W]$ D--compact.
Let $P_+:=1_{[0,\infty)}(D)$. Then the Toeplitz operator 
$P_+WP_+:\im P_+\to \im P_+$ is Fredholm and
\[ 
    \ind(P_+WP_+)=\SF\Bigl((1-s)D+sWDW^*\Bigr)_{0\le s\le 1}.
\]
\end{prop}
\begin{remark}
Note that $P_+WP_+$ is Fredholm on $\im P_+$ if and only
if 
\begin{equation}
P_+WP_++(I-P_+)=I+(W-I)P_+-[W,P_+]P_+
\end{equation}
is Fredholm on $H$ 
with the same index. 
Since $[W,P_+]$ is compact we conclude
\begin{equation}
      \ind(P_+W P_+)=\ind(I+(W-I)P_+)
\end{equation}
which gives \cite[Prop. 3.1]{Bun:SFF}.

To see the compactness of $[W,P_+]$ let
$P_+(WDW^*):=1_{[0,\infty)}(WDW^*)=WP_+W^*$. Then
\begin{equation}
     [W,P_+]W^*=WP_+W^*-P_+=P_+(WDW^*)-P_+.
\end{equation}
Since
$WDW^*-D=[W,D]W^*$ is $D$--compact the operator
$P_+(WDW^*)-P_+$ is compact in view of Corollary \plref{S-compactness}.
\end{remark}

\begin{proof}[Proof of Proposition \plref{Toeplitz}]  
   By assumption we have
\begin{equation}
    D_s:=(1-s)D+sWDW^*=D+s[W,D]W^*.
\end{equation}
Hence we may apply Theorem \plref{S-SFpair} and find that
$(WP_+W^*,P_+)$ is a Fredholm pair and 
\begin{equation}\begin{split}
   \SF(D_s)&=\ind(WP_+W^*,P_+)\\
           &=\ind(P_+:\im WP_+W^*\to\im P_+)\\
           &=\ind(P_+WP_+:\im P_+\to \im P_+)
                \end{split}
\end{equation}
and the result is proved.
\end{proof}

\section{The Theorem of Cordes--Labrousse revisited}\label{secfour}

\begin{prop}\label{S-4.1} $\cB^\sa$ is open and dense in 
$(\cfsa,d_G)$ and also in $(\cfsa,d_R)$. 
Moreover, the topology induced by the graph resp. Riesz metric on $\bsa$ coincides with
the norm topology.
\end{prop}

That the graph metric induces the norm topology on bounded operators is due
to Cordes and Labrousse \cite[Addendum]{CorLab:IIM} who also observed
that the bounded operators are open in the graph metric.
That $\cB^\sa$ is dense in $(\cfsa,d_G)$ was observed in
\cite[Prop. 1.6]{BooLesPhi:UFO}. 

\begin{proof} By Proposition  \plref{S1.1} we know that the Riesz metric is stronger than
the graph metric and applying Proposition  \plref{S1.1} with $D=0$
(which is not excluded!)  we see that the natural inclusion
$\bsa\hookrightarrow (\cfsa,d_R)$ is continuous.  Hence it
suffices to show that $\bsa$ is open in the graph topology, that
$\bsa$ is dense in $(\cfsa,d_R)$ and that the topology induced
by the graph metric on $\bsa$ coincides with the norm topology.

\enum{1}\ Fix $R>0$ and $T\in\bsa$, $\|T\|\le R$. Consider
$\widetilde T\in\csa$ with $d_G(T,\widetilde T)<\frac 12 (1+R)^{-1}$. Then
\begin{equation}
   (\widetilde T+i)^{-1}=(T+i)^{-1}\Bigl(I-(T+i)\bigl((T+i)^{-1}-(\widetilde
   T+i)^{-1}\bigr)\Bigr)
\end{equation}
is invertible with bounded inverse since
\begin{equation}
    \Bigl\|(T+i)\Bigl( (T+i)^{-1}-(\widetilde T+i)^{-1}\Bigr)\Bigr\|\le (1+R)d_G(T,\widetilde
    T)<\frac 12,
\end{equation}
and
\begin{equation}
    (\widetilde T+i)=\sum_{n=0}^\infty\Bigl((T+i)\bigl((T+i)^{-1}
       -(\widetilde T+i)^{-1}\bigr)\Bigr)^n(T+i).\label{G2.3}
\end{equation}
Hence the ball $B_{d_G}(T,\frac 12(1+R)^{-1})$ is contained in
$\bsa$ which proves that $\bsa$ is open in $(\csa,d_G)$.

To prove that it is dense even with respect to the Riesz metric, 
we consider $T\in\csa$ and denote by
$(E_\gl)_{\gl\in\R}$ the spectral resolution of $T$. Let $f_n$ be the
function sketched in Figure \ref{figureone}.

%%%%%%%%%%%%%%%%%%%%%%%%%%%%%%%%%%%%%%%%%%%%%%%%%%%%%%%%%%%%%%%%%%%%%
% Figure 1
%%%%%%%%%%%%%%%%%%%%%%%%%%%%%%%%%%%%%%%%%%%%%%%%%%%%%%%%%%%%%%%%%%%%%
\setlength{\unitlength}{1.0cm}
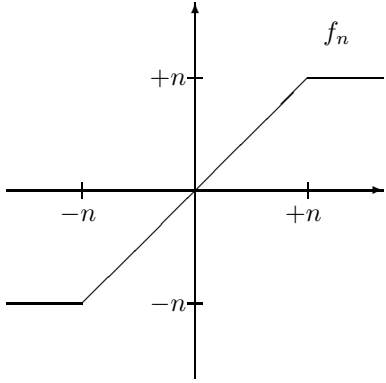
\begin{figure}
\begin{picture}(5.0,5.0)
\put(2.5,0){\vector(0,1){5.0}}   %Koordinatenkreuz
\put(0,2.5){\vector(1,0){5.0}}
\put(0,1){\line(1,0){1}}     %Funktionsverlauf
\put(1,1){\line(1,1){3}}
\put(4,4){\line(1,0){1}}
\put(2.4,1){\line(1,0){0.2}} %Koordinatenmarkierung
\put(1.9,0.9){\parbox{5mm}{$-n$}}  %Koordinatenbeschriftung
\put(2.4,4){\line(1,0){0.2}}
\put(1.9,3.9){\parbox{5mm}{$+n$}}
\put(4.2,4.5){\parbox{5mm}{$f_n$}}
\put(1,2.4){\line(0,1){0.2}}
\put(0.7,2.1){\parbox{5mm}{$-n$}}
\put(4,2.4){\line(0,1){0.2}}
\put(3.7,2.1){\parbox{5mm}{$+n$}}
\end{picture}
\caption[Function $f_n$]{
\label{figureone} Function $f_n$. $f_n(T)$ is bounded and converges
in the graph distance to the operator $T$.}
\end{figure}

We put
\begin{equation}
     f_n(T)=\int_{[-n,n]} \gl dE_\gl+\int_{|\gl|>n}n (\sgn\gl) dE_\gl\in\bsa
\end{equation}
%and find
%\begin{equation}\begin{split}
%    \|(&T+i)^{-1}-(f_n(T)+i)^{-1}
%    \|\\
%   &=\bigl\|\int_{|\gl|>n}(\gl+i)^{-1}-(n(\sgn\gl)+i)^{-1}dE_{\gl}\bigr\|\le \frac{2}{n+1},\end{split}
%\end{equation}
and find
\begin{equation}\begin{split}
    \|F(&T)-F(f_n(T))\|\\
   &=\bigl\|\int_{|\gl|>n}\frac{\gl}{\sqrt{1+\gl^2}}-\frac{n(\sgn\gl)}{\sqrt{1+n^2}}dE_{\gl}\bigr\|\\
   &\le \sup\limits_{|\gl|>n}\bigl|
     \frac{\gl}{\sqrt{1+\gl^2}}-\frac{n(\sgn\gl)}{\sqrt{1+n^2}}\bigr|\\
   &\le \bigl| \frac{n}{\sqrt{1+n^2}}-1  \bigr|,\end{split}
\end{equation}
hence $\lim\limits_{n\to\infty}d_R(T,f_n(T))=0$.

\medskip
\enum{2}\ Let $T\in\bsa,\|T\|\le R$. Then, for $\widetilde T\in\bsa$ with
$d_G(T,\widetilde T)<\frac 12 (1+R)^{-1}$ we have in view of \eqref{G2.3}
\begin{equation}
    \|T-\widetilde T\|\le \sum_{n=1}^\infty (1+R)^{n+1}d_G(T,\widetilde T)^n
                 \le 2(1+R)^2d_G(T,\widetilde T).\label{G2.6}
\end{equation}
Conversely, if $\|T-\widetilde T\|<\frac 12$ we find
\begin{equation}  \begin{split}
           (\widetilde T+i)^{-1}&=\bigl(I-(T+i)^{-1}(T-\widetilde
           T)\bigr)^{-1}(T+i)^{-1}\\
                &=\sum_{n=0}^\infty\bigl((T+i)^{-1}(T-\widetilde T)\bigr)^n
                      (T+i)^{-1}
          \end{split}\label{G2.7}
\end{equation}
and hence
\begin{equation}
     d_G(T,\widetilde T)\le\sum_{n=1}^\infty \|T-\widetilde T\|^n\le 2\|T-\widetilde
     T\|.\label{G2.8}
\end{equation}
\eqref{G2.6} and \eqref{G2.8} show that the topologies induced by $d_G$
and $\|\cdot\|$ on $\bsa$ coincide.
\end{proof}

\eqref{G2.6} and \eqref{G2.8} show a bit more. Namely, given
$T_0\in\bsa$, $R:=\|T_0\|+1$ put $r:=\frac 14 (1+R)^{-2}$. Then
the ball $B_{d_G}(T_0,r)$ is open in the graph and the norm topology. For
$T\in B_{d_G}(T_0,r)$ we find in view of \eqref{G2.6}
\begin{equation}
  \|T-T_0\|\le 2(1+R)^2r=\frac 12
\end{equation}
thus $\|T\|\le R$. Hence \eqref{G2.6} and \eqref{G2.8} may be applied to
arbitrary $T,\widetilde T\in B_{d_G}(T_0,r)$ and we find
\begin{equation}
   \frac 12 d_G(T,\widetilde T)\le \|T-\widetilde T\|\le 2 (1+R)^2d_G(T,\widetilde
   T).
\end{equation}
Hence the norm distance and the $d_G$--distance are equivalent
on $B_{d_G}(T_0,r)$.

Still, the norm distance and the $d_G$--distance are not
globally equivalent! The reason is that $\bsa$ is norm
complete and at the same time $d_G$--dense in $\csa$.

\subsection{The stability of the index}

We are going to present a concise proof of the Theorem of Cordes--Labrousse
on the stability of the index in a very general context.

Let $H:=H^+\oplus H^-$ be a $\Z_2$--graded Hilbert space with grading
operator
\begin{equation}
    \ga=\begin{pmatrix} 1&0\\ 0 &-1\end{pmatrix}.
\end{equation}
We treat the $G$--equivariant index and the Clifford index simultaneously:

\subsubsection*{Case I} Let $G$ be a compact Lie group and let
$\varrho:G\to \cU^{\ev}$ be a unitary representation of $G$ into the
space of even operators on $H$, i.e. $\varrho(g)\alpha=\alpha\varrho(g)$
for $g\in G$.

Spaces of \emph{odd} $G$--equivariant operators are denoted by
a subscript $G$, e.g. $\cfsa_G$, etc. Here, an operator $T$
is called $G$--equivariant if it commutes with $\varrho(g), g\in G$.

\newcommand{\Cl}{C\ell}
\subsubsection*{Case II} Denote by $\Cl_n$ the real Clifford algebra
(Lawson and Michelsohn \cite[Chap. I]{LawMic:SG}), 
i.e. $\Cl_n$ is the universal real $C^*$--algebra
generated by unitaries $e_1,\ldots, e_n$ subject to the relations
\begin{equation}
       e_ie_j+e_je_i=-2\delta_{ij}.
\end{equation}
$\Cl_n$ is $\Z_2$--graded with the generators $e_j$ being of odd degree.

Let $\varrho:\Cl_n\to \cB$ be a faithful unital graded $*$--representation
of $\Cl_n$ on $\cB$.

Spaces of \emph{odd} $\Cl_n$--invariant operators are denoted by a
subscript $n$, for example $\cfsa_n$, etc.

\bigskip
In both cases we now consider an \emph{odd} $\varrho$--equivariant
self--adjoint Fredholm operator $T\in \cfsa_G$ ($\cfsa_n$).
Then $\ker T$ is a $\Z_2$--graded $\varrho$--module.

\subsubsection*{Case I} Denote by $(\ker T)^\pm:=\ker T\cap \ker (\alpha\mp I)$
the $\pm$--part of $\ker T$. Then $(\ker T)^\pm$ are $G$--modules. One
puts
\begin{equation}
     \ind_G(T):=[\ker T]:=[(\ker T)^+]-[(\ker T)^-]\in R(G),
\end{equation}
where $R(G)$ is the ring of virtual finite--dimensional representations
of $G$, i.e. it is the Grothendieck group of the semiring of equivalence
classes of finite--dimensional representations.

\subsubsection*{Case II} Following Atiyah, Bott, and Shapiro
\cite{AtiBotSha:CM}
(cf. also \cite[Sec. I.9]{LawMic:SG}) 
let $\hat\cM_n$ be the Grothendieck group
of equivalence classes of finite--dimensional $\Z_2$--graded $\Cl_n$--modules.
Then there is a canonical isomorphism
\begin{equation}
     \hat \cM_n/\hat\cM_{n+1}\simeq KO^{-n}(pt).
\end{equation}
Recall that $KO^{-n}(pt)$ is $8$--periodic, $KO^0(pt)\simeq KO^{-4}(pt)\simeq \Z,
KO^{-1}(pt)\simeq KO^{-2}(pt)\simeq \Z_2, $ and the remaining groups
vanish.

The isomorphism
\begin{equation}
    \hat \cM_0/\hat \cM_1\simeq KO^0(pt)\simeq\Z \label{G2.15}
\end{equation}
is given by sending the graded vector space $V$ (a $\Cl_0$--module)
to its graded dimension $\dim_{\Z_2} V:=\dim V^+-\dim V^-$.

Again, $\ker T$ is a $\Z_2$--graded $\Cl_n$--module and one puts
\begin{equation}
   \ind_n T:= [\ker T]\in \hat\cM_n/\hat\cM_{n+1}.
\end{equation}
Note that an odd self--adjoint Fredholm operator $T$ takes the form
\begin{equation}T=
\begin{pmatrix}0& (T^+)^*\\ T^+ &0\end{pmatrix}
\end{equation}
and in view of \eqref{G2.15} $\ind_0 T$ is nothing but the ordinary Fredholm
index of $T^+$.

\begin{lemma}\label{S2.2} Let $T\in\cfsa_G$ or $T\in\cfsa_n$. Let $\eps>0$
such that $\pm\eps\not\in\spec T$ and $[-\eps,\eps]\cap\specess(T)=\emptyset$.
Then $\im 1_{[-\eps,\eps]}(T)$ is a $G$--module ($\Cl_n$--module)
and
\[ \bigl[\im(1_{[-\eps,\eps]}(T))\bigr]=\bigl[\ker T\bigr]
\]
in $R(G)$ resp. $\hat \cM_n/\hat \cM_{n+1}$.
\end{lemma}
\begin{proof} We first note that the choice of $\eps>0$ is possible. 
Namely, since $T$ is a Fredholm operator $0$ is
not in the essential spectrum. Hence,
$0$ is at most an isolated point of $\spec T$. 
Since spectrum and essential spectrum
are closed one may choose $\eps>0$ as stated.

Abbreviate $V:= \im(1_{[-\eps,\eps]}(T))$. Then we have a
$\varrho$--equivariant decomposition
\begin{equation}
    V=\ker T\bigoplus\mathop{\oplus}_{0<\gl<\eps}\ker(T^2-\gl^2).
   \label{G2.19}
\end{equation}
Now consider $\gl>0$:
\subsubsection*{Case I} $|T|^{-1}T:\ker(T^2-\gl^2)^\pm\to
\ker(T^2-\gl^2)^\mp$ is a $G$--equivariant isomorphism, hence
$[\ker(T^2-\gl^2)]=0$ in $R(G)$.

\subsubsection*{Case II} On $\ker(T^2-\gl^2)$ consider $J:=|T|^{-1}T$.
$J$ is odd, unitary, and $J^2=I$. Thus with respect to the grading
it takes the form
\begin{equation}\begin{pmatrix}
      0 & (J^+)^*\\ J^+&0\end{pmatrix}.
\end{equation}
Put
\begin{equation} E_{n+1}:=\begin{pmatrix}
      0 & (-J^+)^*\\ J^+&0\end{pmatrix}.
\end{equation}
Again, $E_{n+1}$ is odd,  unitary, and $E_{n+1}^2=-I$. Moreover, $E_{n+1}$
anticommutes with $\varrho(e_k), k=1,...,n$. Hence
$\varrho(e_1),...,\varrho(e_n),E_{n+1}$ make $\ker(T^2-\gl^2)$ into a
graded $\Cl_{n+1}$--module and hence $[\ker(T^2-\gl^2)]=0$ in
$\hat \cM_n/\hat \cM_{n+1}$.

In view of \eqref{G2.19} we thus have in both cases
$[\im(1_{[-\eps,\eps]}(T))]=[\ker T]$.
\end{proof}

For the ordinary Fredholm index 
the following Theorem is due to Cordes and Labrousse
\cite{CorLab:IIM}. For bounded operators in the current equivariant
context it can be found in \cite[Sec. III.10]{LawMic:SG}. 
However, since the Riesz
topology and the graph topology are different (Proposition  \plref{S1.2}), 
(10.8) in \cite{LawMic:SG} 
is problematic and valid only in the context
of unbounded operators with a fixed domain.

\begin{theorem}\label{stability} The $G$--index
\[ \ind_G:\cfsa_G\to R(G)\]
and the Clifford index
\[ \ind_n:\cfsa_n\to KO^{-n}(pt)\]
are locally constant with respect to the graph topology on $\cfsa$.
\end{theorem}
\begin{proof} Fix a $T\in \cfsa_G$ (resp. $T\in\cfsa_n$). Since
$T$ is a Fredholm operator, $0\not\in\specess(T)$. Thus there is
an $\eps>0$ such that $\pm\eps\not\in\spec T$ and
$[-\eps,\eps]\cap\specess(T)=\emptyset$.

Moreover, in view of \cite[Prop. 2.10]{BooLesPhi:UFO} there is an
open neighborhood $\cN\subset\cfsa_G$ ($\cfsa_n$) of $T$
such that $\cN\ni S\mapsto 1_{[-\eps,\eps]}(S)$ is continuous and
finite--rank projection valued.

In particular, making $\cN$ smaller if necessary, we may assume that
\begin{equation}
    \| 1_{[-\eps,\eps]}(S)- 1_{[-\eps,\eps]}(T)\|<1
\end{equation}
for all $S\in\cN$.

It is well--known that if two orthogonal projections $P,Q$ satisfy
$\|P-Q\|<1$ then $P$ maps $\im Q$ isomorphically onto $\im P$.
Hence $1_{[-\eps,\eps]}(S)$ is a $\varrho$--equivariant isomorphism
from $\im(1_{[-\eps,\eps]}(T))$ onto $\im(1_{[-\eps,\eps]}(S))$.
Thus by Lemma \plref{S2.2}
\begin{equation}
  [\ker S]=[\im(1_{[-\eps,\eps]}(S))]=[\im(1_{[-\eps,\eps]}(T))]=[\ker T]
\end{equation}
for $S\in\cN$ and the Theorem is proved.
\end{proof}

\section{Uniqueness of the spectral flow}\label{secfive}

\subsection{The general set--up} \label{secfiveone}
We start fixing some basic notation und introducing
the problem: 

\begin{dfn}\label{defhomotopy} For a topological space $X$
and a subspace $Y\subset X$ we denote by $\Omega(X,Y)$ the set
of paths $f:[0,1]\to X$ with endpoints in $Y$.
Instead of $\Omega(X,X)$ we also write $\Omega(X)$. Paths
are always assumed to be continuous.
 
Paths $f,g\in\Omega(X,Y)$ are (free) homotopic if there is
a continuous map $H:[0,1]\times [0,1]\to X$ with the properties
\begin{enumerate}
\item $H(0,.)=f, H(1,.)=g$,
\item $H(s,0)\in Y, H(s,1)\in Y$ for all $s$.
\end{enumerate}
The set of homotopy classes in this sense is
denoted by $\pitilde_1(X,Y)$.
\end{dfn}

Note that we do not require a base point in $Y$ to stay fixed during
the deformation, however endpoints are only allowed to move within $Y$.
Therefore $\pitilde_1(X,Y)$ is not
the relative homotopy set usually introduced in algebraic topology
textbooks. If $y_0$ is a base point in $Y$ 
then the relative homotopy set is denoted
by $\pi_1(X,Y,y_0)$.

\begin{dfn}\label{defproblem} Let $X$ be a topological space and $Y\subset X$
a subspace. For a map
\begin{equation}
    \mu:\Omega(X,Y)\longrightarrow \Z
\end{equation}
the properties \emph{Concatenation} and \emph{Homotopy}
are defined as follows:
\begin{enumerate}
\item\emph{Concatenation:} If $f,g\in \Omega(X,Y)$
are paths with $f(1)=g(0)$ then
\begin{equation}
    \mu(f*g)=\mu(f)+\mu(g).
\end{equation}
\item\emph{Homotopy:} $\mu$ descends to a map
$\mu:\pitilde_1(X,Y)\to\Z$.\label{homotopyfd}
\end{enumerate}
\end{dfn}

\begin{lemma}\label{constantcurve} Let $X,Y$ be as before and
let $ \mu:\Omega(X,Y)\longrightarrow \Z$ be a map which satisfies
\emph{Concatenation} and \emph{Homotopy}. 

Then for each $x_0\in Y$ the restriction of $\mu$ to $\Omega(X,x_0)$
is a homomorphism $\pi_1(X,x_0)\to \Z$.
Furthermore, for any path $f:[0,1]\to Y$ we have $\mu(f)=0$.
\end{lemma}
\begin{proof} The first claim is obvious. To prove the second 
claim we first note that a constant path $f=f(0)$ satisfies 
$\mu(f)=\mu(f*f)=\mu(f)+\mu(f)$, thus $\mu(f)=0$.

A general path $f\in\Omega(Y)$ is homotopic in $\Omega(X,Y)$ to
the constant path $f(0)$ via the homotopy $f_s(t):=f(st)$ and 
we reach the conclusion.
\end{proof}

\subsection{Warm--up: Uniqueness in the bounded case} \label{secfivetwo}
Let $H$ be a separable complex Hilbert space, i.e. the Hilbert
dimension is finite or countably infinite. If $A\subset \cC$ is
a set of operators in $H$ we denote by $GA$ the set of
invertible elements in $A$ (with bounded inverse), by
$\cF A$ the Fredholm operators in $A$, and by $\cF_* A$ the Fredholm operators
in $A$ which are neither essentially positive nor essentially negative.

The following uniqueness--theorem for 
the spectral flow in the classical situation of bounded operators
follows easily from the isomorphism \eqref{G-I.3}. This is of course
folklore.

\begin{theorem}\label{uniqueness-bounded} Let $H$ be an infinite--dimensional
separable complex Hilbert space and let 
$\mu:\Omega(\bfsa_*,G\bfsa_*)\longrightarrow \Z$ be a map which satisfies
\emph{Concatenation}, \emph{Homotopy} and
\begin{enumerate}
\item[\relax] \emph{Normalization:} There is a $T_0\in G\bfsa_*$ and
a rank one orthogonal projection
$P\in\bsa$ commuting with $T_0$ such that the path 
\[f_P(t):=tP+(I-P)T_0,\quad -1/2\le t\le 1/2\]
satisfies
\[  \mu(f_P)=1.
\]%\label{normalizationbd}
\end{enumerate}
Then $\mu$ equals the spectral flow.  
\end{theorem}
\begin{proof} 
\comment{We first note that it is a consequence of Kuiper's Theorem
\cite{Kui:HTU} that $G\bfsa_*$ is contractible. Let us briefly sketch this:
by pushing the spectrum to $\pm 1$ we see that $G\bfsa_*$ is homotopy
equivalent to the space $\bigsetdef{T\bsa}{T^2=1, \dim\ker(T\pm I)=\infty}$.
On this space the unitary group $\cU$ acts transitively by conjugation
with stabilizer isomorphic to $\cU(\ker (T_0-I)\times \cU(ker(T_0+I))$ for any
fixed $T_0$. By Kuiper's Theorem $\cU, \cU(\ker (T_0-I),$ and $\cU(ker(T_0+I)$
are contractible and hence $G\bfsa_*$ must be contractible, too.

In particular the fundamental group of $G\bfsa_*$ vanishes and the long
exact homotopy sequence of the pair $(\bfsa_, G\bfsa_*)$ shows that
for any fixed $T_0\in G\bfsa_*$ the natural map
}
We first note that the spectral flow satisfies \emph{Concatenation, Homotopy}, and
\emph{Normalization}.

$G\bfsa_*$ is connected. Therefore, we may choose
a path $g_P\in\Omega(G\bfsa_*)$ from $f_P(1)$ to $f_P(0)$. In view of
Lemma \plref{constantcurve}, and \emph{Normalization} 
we thus have $\mu(g_P*f_P)=\SF(g_P*f_P)=1$. By \eqref{G-I.3} the closed
path $g_P*f_P$ must be a generator of $\pi_1(\bfsa_*,f_P(0))$
and consequently $\mu=\SF$ on $\pi_1(\bfsa_*,f_P(0))$ again by
Lemma \plref{constantcurve}.

If $f\in\Omega(\bfsa_*,G\bfsa_*)$ is arbitrary we choose paths 
$g_1, g_2:[0,1]\to G\bfsa_*$ with $g_1(0)=f_P(0), g_1(1)=f(0), g_2(0)=f(1),
g_2(1)=f_P(0)$. Then the path $g_1*f*g_2$ is closed and Lemma \plref{constantcurve}
yields
\begin{equation}
   \mu(f)=\mu(g_1*f*g_2)=\SF(g_1*f*g_2)=\SF(f)
\end{equation}
and we are done.
\end{proof}
   
Amazingly the finite--dimensional analogue of the previous Theorem is slightly more
complicated due to the fact that in this case $G\bsa$ is not connected. Of course,
if $H$ is finite--dimensional then $\bfsa=\bsa$.

\begin{prop} \label{S3.1fd}
If $\dim H<\infty$ then the path components of $G\bsa$ are
labelled by $\rank(1_{[0,\infty)}(T))\in\{0,\ldots,\dim H\}$.
\end{prop}

\begin{proof} \ $T\mapsto \rank(1_{[0,\infty)}(T))$ is continuous
on $G\bsa\subset GL(\dim H)$ and maps onto $\{0,...,\dim H\}$.

Obviously, for $T\in G\bsa$ there is a path in $G\bsa$
connecting $T$ with $2P-I$ for $P=1_{[0,\infty)}(T)$
which shows injectivity.
\end{proof}

\begin{lemma} \label{S3.2} Let $\dim H<\infty$ and let
$f,g\in\Omega(\bsa,G\bsa)$ be paths with the same initial points
$f(0)=g(0)$.

Then $f,g$ define the same class in $\pitilde_1(\bsa,G\bsa)$ if and
only if
\begin{equation}
\rank\bigl(1_{[0,\infty)}(f(1))\bigr)=\rank\bigl(1_{[0,\infty)}(g(1))\bigr).
\label{G3.1}
       \end{equation}
\end{lemma}
\begin{proof} If $f,g$ define the same class in $\pitilde_1(\bsa,G\bsa)$
then $f(1)$ and $g(1)$ lie certainly in the same path component of
$G\bsa$ and hence \eqref{G3.1} holds
by Proposition  \plref{S3.1fd}.

The exact homotopy sequence of the pair
$(\bsa,G\bsa)$ gives a bijection
\begin{equation}
    \pi_1(\bsa,G\bsa,f(0))\to\pi_0(G\bsa), \quad [h]\mapsto [h(1)],
\end{equation}
hence from Proposition  \plref{S3.1fd} we infer
that $f,g$ even define the same class in the relative homotopy set
$\pi_1(X,Y,f(0))$, in particular they define the same class in $\pitilde_1(X,Y)$.
\end{proof}

\begin{theorem}\label{S3.3}\label{uniqueness-fd} Let $H$ be a finite--dimensional Hilbert space and let
\[\mu:\Omega(\bsa,G\bsa) \longrightarrow \Z\]
be a map which satisfies \emph{Concatenation}, \emph{Homotopy} and
\emph{Normalization} in the following sense:
\begin{enumerate}
\item[\relax] There is a rank one orthogonal projection $P\in\bsa$ such that for all $A\in\bsa$
\begin{equation}
        \mu\bigl((tP+(I-P)A(I-P))_{-1/2\le t\le 1/2}\bigr)=1.
\end{equation}\label{normalizationfd}
\end{enumerate}
Then
\[
\mu(f)
  =\rank(1_{[0,\infty)}(f(1)))-\rank(1_{[0,\infty)}(f(0)))
  =\SF(f)
\]
for all $f\in\Omega(\bsa,G\bsa)$.
\end{theorem}
\begin{proof} First note that \emph{Normalization} holds for any rank one orthogonal projection:
namely all rank one orthogonal projections are unitarily equivalent and
the unitary group is connected, hence \emph{Homotopy} implies that \emph{Normalization} holds for any
rank one orthogonal projection $P$.

Now consider a path $f\in\Omega(\bsa,G\bsa)$. In view of Proposition 
\plref{S3.1fd} and \emph{Homotopy} we may assume that
\begin{equation}   f(0)=2P-I \text{ and }  f(1)=2Q-I,
\end{equation}
where $P,Q$ are orthogonal projections. Put
\begin{equation}
    \gamma_P(t)=2P-I+2t(I-P), \quad 0\le t\le 1.
\end{equation}
Then $\gamma_P^-*f*\gamma_Q$ starts and ends at I.

By Lemma \plref{S3.2} $\gamma_P^-*f*\gamma_Q$ is homotopic in $\Omega(\bsa,G\bsa)$
to the constant curve I, hence from \emph{Concatenation}, \emph{Homotopy} and Lemma
\plref{constantcurve} we infer that 
\begin{equation}
    \mu(f)=\mu(\gamma_P)-\mu(\gamma_Q).
\end{equation}
%It remains to show $\mu(\gamma_P)=n-\rank P$.
If we can show that $\mu(\gamma_P)=n-\rank P$ then
we find 
\begin{equation}
\mu(f)=\rank Q-\rank
P=\rank\bigl(1_{[0,\infty)}(f(1))\bigr)-\rank\bigl(1_{[0,\infty)}(f(0))\bigr)
=\SF(f).\end{equation}
In the last equation we have used Corollary \plref{S-3.7}.

It remains to show $\mu(\gamma_P)=n-\rank P$.
Fix an orthonormal basis such that $P$ has the matrix representation
\begin{equation}
P=\begin{pmatrix} I_k&0\\0&0\end{pmatrix}.\end{equation}
Then $\gamma_P$ is homotopic to $\gamma_1*\gamma_2*\ldots*\gamma_{n-k}$
where
\begin{equation}
   \gamma_j(t)=\begin{pmatrix}I_{k+j-1}&&\mathbf{0}\\
                               &2t-1&\\
                           \mathbf{0}&&-I_{n-k-j}\end{pmatrix}.
\end{equation}
\emph{Normalization} and the remark at the beginning of this proof show
$\mu(\gamma_j)=1$. Hence we find with \emph{Concatenation} and \emph{Homotopy}
\begin{equation}
    \mu(\gamma_P)=\sum_{j=1}^{n-k}\mu(\gamma_j)=n-k=n-\rank P.\qedhere
\end{equation}
\end{proof}

\subsection{Uniqueness of the spectral flow for graph continuous paths}
\label{secfivethree}

Let $H$ be an infinite--dimen\-sional
separable complex Hilbert space. During this subsection we consider the \emph{graph
topology} on $\csa$.

We treat the bounded and the unbounded case simultaneously. Thereby we reprove
Theorem \plref{uniqueness-bounded} without using \eqref{G-I.3}.
We now let X be 
\[
\bfsa_*:=\bigsetdef{T\in\cB}{T=T^*, T\text{ Fredholm}, \specess T\cap
\R_\pm\not=\emptyset} \text{ or } \cfsa.
\]
$X$ is connected since $\bfsa_*$ is trivially connected and $(\cfsa,d_G)$
is connected by \cite[Thm. 1.10]{BooLesPhi:UFO}.

Next put $Y:=GX$, i.e. $Y=G\bfsa_*=\bfsa_*\cap G\bsa$ or $Y=G\csa$. Again,
$Y$ is connected. In the bounded case this is trivial. In the 
unbounded case it is less obvious:

\begin{prop}\label{S3.1id}
If $\dim H=\infty$ then $(G\csa,d_G)$ is path connected.
\end{prop}

This is proved in the spirit of \cite[Thm. 1.10]{BooLesPhi:UFO} and in
fact we also reprove the connectedness of $\cfsa$ in a slightly
different way. During the proof we use the notation of loc. cit. freely.
\begin{proof}
We look at the Cayley picture and consider
$U=\kappa(T)$. Recall that the Cayley transform $\kappa$ is a
homeomorphism from $G\csa$ onto
\begin{equation}
\kappa(G\csa)=\bigsetdef{U\in\cU}{U+I
\text{ invertible and } U-I \text{ injective }}.
\end{equation}

As in loc. cit. $H=H_+\oplus H_-$ is the direct sum of the
spectral subspaces of $U$ corresponding to $\bigsetdef{\gl\in
S^1}{\Im \gl\ge 0}$ and $\bigsetdef{\gl\in S^1}{\Im \gl< 0}$ and
by squeezing the spectrum down to $+i$ and $-i$ one can deform
$U$ within $\kappa(G\csa)$ to
\begin{equation}
    U_1=+i I_+\oplus -i I_-
\end{equation}
(cf. Figure \ref{figtwo}). 
Now since $\dim H=\infty$ we have $\dim H_+=\infty$ or $\dim
H_-=\infty$. 
%%%%%%%%%%%%%%%%%%%%%%%%%%%%%%%%%%%%%%%%%%%%%%%%%%%%%%%%%%%%%%%%%%%%%%%%%%%%%%%%
%
%               FIGURE 2
%
% Due to memory restrictions on some tex systems, we provide two
% alternatives for the inclusion of the figure 2
%
% Caption, needed in both alternatives
% 
\newcommand{\captionoffiguretwo}{
\caption[Connecting a fixed $U\in \kappa(G\cfsa)$ to $iI$.]{
Connecting a fixed $U\in\kappa(G\cfsa)$ to $iI$. 
Case II (infinite rank $U_+$) is first deformed to $-i I$ and then
Case I (infinite rank $U_-$) applies. Through the deformation $-1$ is
never a spectral point.}
\label{figtwo}
}% end of figuretwocaption
%
%
%  
%  Alternative 1: Figure is calculated on the fly (needs a lot of memory and a fast machine)
%\input{Lesch-figure2.tex} %contains the figure
%\begin{figure}
%\[
%\figuretwo          % figuretwo is defined in figure2.tex
%\]
%\captionoffiguretwo
%\end{figure}
%
%
%  Alternative 2: Figure is included as eps, needs figure2.eps
\begin{figure}
\[
\epsfig{file=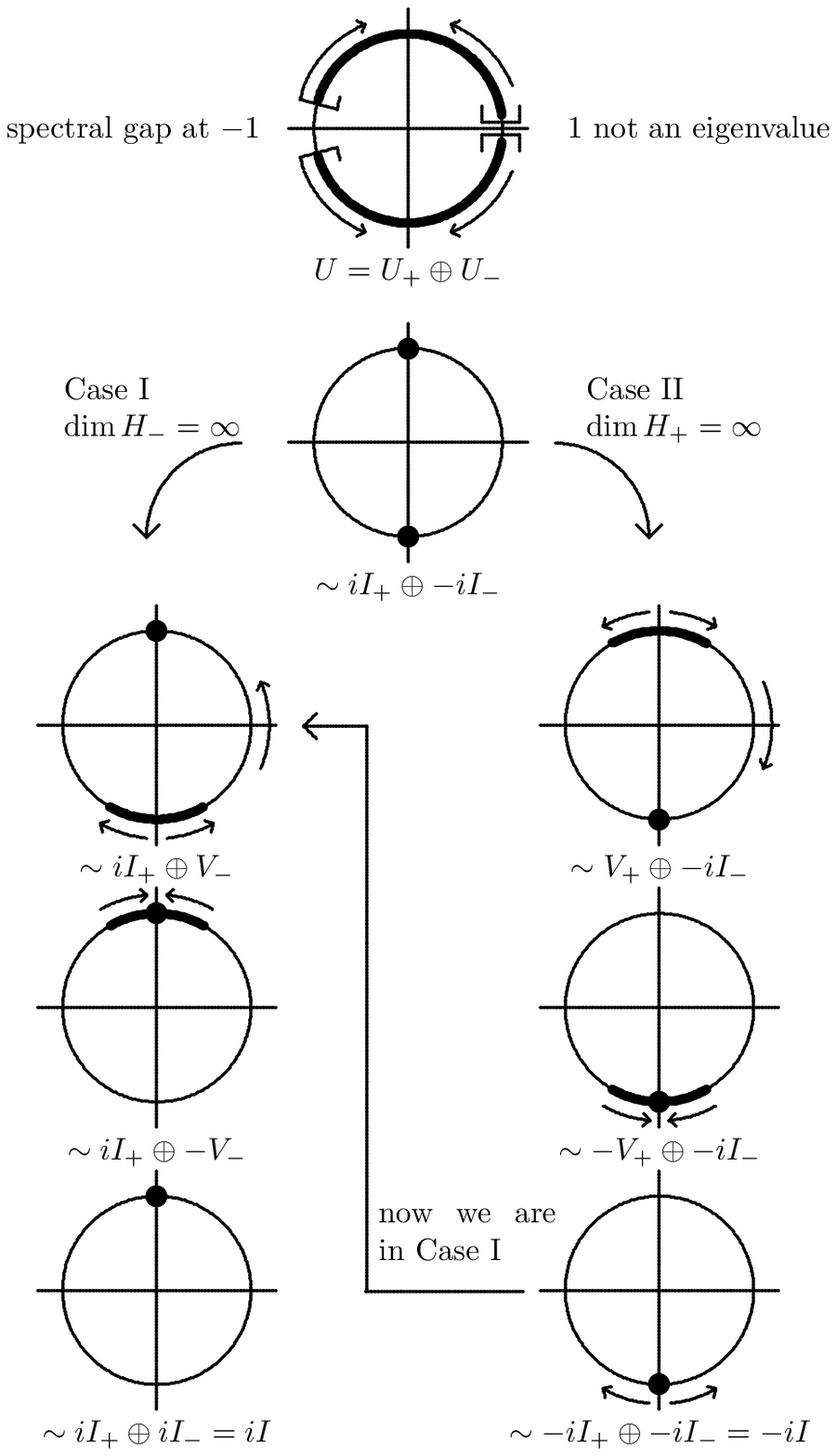}   
% remember to include the package "epsfig" in the beginning of the document
\]
\captionoffiguretwo
\end{figure}
%%%%%%%%%%%%%%%%%%%%%%%%%%%%%%%%%%%%%%%%%%%%%%%%%%%%%%%%%%%%
% end of inclusion of Figure 2!!
%%%%%%%%%%%%%%%%%%%%%%%%%%%%%%%%%%%%%%%%%%%%%%%%%%%%%%%%%%%%%

\subsubsection*{Case I: $\dim H_-=\infty$}
As described in loc. cit. we may un-contract $-iI_-$ in such a way that no eigenvalues
remain, i.e.
\begin{equation}
     U_1\sim iI_+\oplus V_-,
\end{equation}
where $\spec V_-$ consists of a little arc centred on $-i$ and
$V_-$ has no eigenvalues. We then rotate this arc up through $+1$
until it is centered on $+i$. Then we contract the spectrum
to be $+i$. 
This homotopy will stay within $\kappa(G\csa)$ and deform $U_1$
to $i I_H$.

\subsubsection*{Case II: $\dim H_+=\infty$} As in Case I we now
un-contract $+i I_+$ and deform $U_1$ into $U_2=-i I_H$. Now the
operator $U_2$ has $\dim H_-(U_2)=\infty$. Applying Case I we deform
$U_2$ to $+i I_H$.
\end{proof}

Next we choose a base point $T_0\in G\bfsa_*\subset Y$ with $\spec T_0=\{\pm 1\}$, 
in particular $T_0$ is bounded.
With these preparations the uniquess of the spectral flow on $\cfsa$
reads as follows:

\begin{theorem}\label{S-3.4}\label{uniqueness-graph} Let
\[\mu:\Omega(X,Y)\longrightarrow \Z\]
be a map which satisfies
\emph{Concatenation}, \emph{Homotopy} and \emph{Normalization} in the following sense:
\begin{enumerate}
\item[\relax] There is a rank one orthogonal projection
$P\in\bsa$ 
such that the operator $(I-P)T_0(I-P)\in\bsa(\ker P)$ is invertible
and such that 
\begin{equation}
        \mu\bigl((tP+(I-P)T_0(I-P))_{-1/2\le t\le 1/2}\bigr)=1.
\end{equation}
\end{enumerate}
Then $\mu$ equals the spectral flow.
\end{theorem}
\begin{proof}
We will deform a general path $f\in\Omega(X,Y)$ in
several steps into some normal form and then compare $\mu$ and $\SF$
on this normal form. The latter problem is basically reduced to the
finite--dimensional case which was treated as a warm--up in the previous
subsection.

In the sequel $\sim$ denotes homotopy in $\Omega(X,Y)$.
Consider $f\in\Omega(X,Y)$.

\newtheorem{assert}{Assertion}
\begin{assert}\label{assertone} There
exist paths $f_1,\ldots,f_n\in\Omega(X,Y)$ having the following
properties
\begin{enumerate}
\item $f\sim f_1*\ldots *f_n$.
\item There exist $\eps_j>0$ such that for all $t\in [0,1]$
we have $\pm \eps_j\not\in\spec f_j(t)$ and $\specess(f_j(t))\cap
[-\eps_j,\eps_j]=\emptyset$.
\end{enumerate}
\end{assert}

This is a basic fact about paths of Fredholm operators and has 
nothing to do with our assumptions on $\mu$. Namely, for each $t$
the operator $f(t)$ is Fredholm and hence $0\not\in\specess f(t)$.
Hence by compactness (cf. \cite[Prop. 2.10 and Def. 2.12]{BooLesPhi:UFO})
there is a subdivision $0=t_0< t_1<\ldots < t_n=1$ of the interval $[0,1]$
and positive real numbers $\eps_j, j=1,\ldots, n$, such that
$\pm\eps_j\not\in\spec f(t)$ and $[-\eps_j,\eps_j]\cap \specess(f(t))
=\emptyset$ for $t_{j-1}\le t\le t_j, j=1,\ldots,n$.

Hence the candidates for the $f_j$ are $f|[t_{j-1},t_j]$ (reparametrized
over [0,1]). The problem is that $f(t_j)$ need not be invertible and
hence need not be in $Y$. However, $0$
is at most an isolated point of the spectrum of $f(t_j)$
since $f(t_j)$ is Fredholm. Hence
there is a $\delta>0$ such that $f(t_j)+s\delta$ is invertible
for $0<s\le\delta, j=1,\ldots,n-1$ and $0\le s\le\delta, j=0,n$.
By compactness we may choose $\delta$ so small that additionally
$\pm\eps_j\not\in\spec(f(t)+s\delta)$ and $[-\eps_j,\eps_j]\cap
\specess(f(t)+s\delta)=\emptyset$ for all $t_{j-1}\le t\le t_j,
0\le s\le \delta$. Thus 
\begin{equation}
     H(s,t):=f(t)+sI,\quad 0\le s\le\delta
\end{equation}
is a homotopy in $\Omega(X,Y)$ and we reach the conclusion with
$f_j:=f|[t_{j-1},t_j]+\delta I$ (reparametrized over $[0,1]$).

\begin{assert}\label{asserttwo} Suppose that there is an $\eps>0$ such
that $\pm\eps\not\in\spec(f(t))$ and $[-\eps,\eps]\cap\specess(f(t))=\emptyset$
for all $t\in[0,1]$. Then $f$ 
is homotopic in $\Omega(X,Y)$ to a path $g$
having the following properties: there is a finite rank orthogonal projection
$\widetilde Q$ and an operator $S\in G\bigl((I-\widetilde Q)X(I-\widetilde Q)\bigr)$
such that with respect to the orthogonal decomposition
$H=\im \widetilde Q\oplus \ker \widetilde Q$
we have
\begin{equation}
      g(t)=\begin{pmatrix} g_0(t) &0\\ 0 & S
             \end{pmatrix}.\label{G3.16}
\end{equation}
\end{assert}

Since $[-\eps,\eps]\cap\specess(f(t))=\emptyset$ the spectral
projection
\begin{equation}
      E(t):=1_{(-\eps,\eps)}(f(t))
\end{equation}
is of finite rank and
$E(t)$ depends continuously on $t$ by \cite[Prop. 2.10]{BooLesPhi:UFO}.

By \cite[Prop. 4.3.3]{Bla:KTO}
there exists a continuous family of unitaries
$U:[0,1]\to \cU, U(0)=I$ such that $E(t)=U(t)E(0)U(t)^*$.

Now consider the homotopy
\begin{equation}
     H(s,t):=U(st)^*f(t)U(st),\quad 0\le s,t\le 1.
\end{equation}
$H$ is certainly a homotopy in $\Omega(X,Y)$. Furthermore,
\begin{equation}\begin{split}
   H(1,t)E(0)&=U(t)^*f(t)U(t)E(0)=U(t)^*f(t)E(t)U(t)\\
             &=U(t)^*E(t)f(t)U(t)=E(0)H(1,t),\end{split}
\end{equation}
since $E(t)$ is a spectral projection of $f(t)$.
Thus the orthogonal projection $E(0)$ 
commutes with the self--adjoint operator $H(1,t)$ and hence 
with respect to the decomposition $H=\im E(0)\oplus \ker E(0)$
the operator $H(1,t)$ takes the form
\begin{equation}
       H(1,t)=\begin{pmatrix} g_0(t)&0\\ 0& g_1(t)\end{pmatrix}.
\end{equation}
By construction $g_1(t)$ is invertible for all $t$ and thus the map
\begin{equation}
       \begin{pmatrix} g_0(t)&0\\ 0& g_1(st)\end{pmatrix},
\quad 0\le s,t\le 1
\end{equation}
homotops $H(1,\cdot)$ to
\begin{equation}
       \begin{pmatrix} g_0(t)&0\\ 0& g_1(0)\end{pmatrix},
\quad 0\le s,t\le 1,
\end{equation}
and Assertion \ref{asserttwo} is proved with $S=g_1(0)$ and $\widetilde Q=E(0)$.

\begin{assert}\label{assertthree}
Consider $g$ as in Assertion \plref{asserttwo}.
Then $g$ is homotopic in $\Omega(X,Y)$ to a path
$h$ having the following properties:
there is a finite rank orthogonal projection $Q\ge P$ such
that $(I-Q)T_0(I-Q)\in\bsa(\ker Q)$ is invertible and such
that with respect to the orthogonal decomposition
$H=\im Q\oplus \ker Q$ we have
\begin{equation}
      h(t)=\begin{pmatrix} h_0(t) &0\\ 0 & (I-Q)T_0(I-Q)
             \end{pmatrix}.\label{G3.20}
\end{equation}
\end{assert}

There is a unitary operator $U\in\cU$ such that
$Q:=U^*\widetilde Q U\ge P$ and $(I-Q)T_0(I-Q)\in\bsa(\ker Q)$ is
invertible. Since the unitary group $\cU$ is
connected we may choose a path $U:[0,1]\to\cU$ with $U(1)=U, U(0)=I$.
Then $H(s,t):=U(s)^*g(t)U(s)$ homotops $g$ in $\Omega(X,Y)$
to $U^*gU$. W.r.t. the orthogonal decomposition
$H=\im Q\oplus \ker Q$ the latter takes the form
\begin{equation}
      \widetilde h(t)=\begin{pmatrix} h_0(t) &0\\ 0 & \tilde S
             \end{pmatrix}.
\end{equation}
Finally, since $G\bigl((I-Q)X(I-Q)\bigr)$ is path connected\footnote{Note that
$Q$ is of finite rank and hence $(I-Q)X(I-Q)$ is either 
$G\bfsa_*(\ker Q)$ or $G\csa(\ker Q)$. In either case 
$G\bigl((I-Q)X(I-Q)\bigr)$ is path connected since $\ker Q$ is a separable
infinite--dimensional Hilbert space (Proposition \plref{S3.1id}).}
there is a path
$h_1:[0,1]\to G\bigl((I-Q)X(I-Q)\bigr)$ with $h_1(0)=\tilde S$ and $h_1(1)=(I-Q)T_0(I-Q)$
and
\begin{equation}
       \begin{pmatrix} h_0(t)&0\\ 0& h_1(s)\end{pmatrix},
\quad 0\le s,t\le 1
\end{equation}
homotops $\widetilde h$ to the claimed path $h$.

\subsubsection*{\textsc{Finish of proof}}\label{finishproof} In view of \emph{Homotopy}
 and \emph{Concatenation}
and in view of Assertions \ref{assertone}--\ref{assertthree} 
it remains to show that
for $h$ in \eqref{G3.20} we have $\mu(h)=\SF(h)$.

Consider the map
\begin{equation}\begin{split}
    \sigma:&\Omega(\bsa(\im Q),G\bsa(\im Q))\longrightarrow\Z,\\
      & \sigma(f):=\mu(\begin{pmatrix}f&0\\0&(I-Q)T_0(I-Q)\end{pmatrix}).
                \end{split}
\end{equation}
$\sigma$ inherits \emph{Concatenation} and \emph{Homotopy} immediately from $\mu$.
$\sigma$ is also normalized since $P\le Q$.

Thus we may apply Theorem \plref{S3.3} and conclude 
\begin{equation}
\mu(h)=\sigma(h_0)=\rank(1_{[0,\infty)}(h_0))-
\rank(1_{[0,\infty)}(h_0))=\SF(h_0)=\SF(h).\qedhere
\end{equation}
\end{proof}

\subsection{Uniqueness of the spectral flow for Riesz continuous paths}\label{secfive-Riesz}

\label{secfivefour}
Because of the next result all results about the spectral flow of paths of bounded
operators carry over verbatim to \emph{Riesz} continuous paths of unbounded operators.
The drawback is, as mentioned in the introduction, that the Riesz metric is so
strong that it is hard to prove continuity of maps into the space.

\begin{theorem}\label{uniqueness-Riesz} The natural inclusion of the pair
\[j:(\bfsa,G\bsa)\hookrightarrow
\bigl((\cfsa,G\csa),d_R\bigr)\]
is a homotopy equivalence.
\end{theorem}

\begin{proof} The image of the
Riesz map $F(T)=T(I+T^2)^{-1/2}$ was determined in \cite[Prop. 1.5]{BooLesPhi:UFO}, i.e.
\begin{equation}
   F(\csa)=\bigsetdef{S\in\bsa}{\|S\|\le 1 \text{ and } S\pm I \text{ both injective}}=:X.
\end{equation}
$F$ is a homeomorphism of $\csa$ onto $X\subset \bsa$ by definition of the Riesz metric.
From the functional calculus we know that 
$F$ maps the (essential) spectrum of $T$ onto
the (essential) spectrum of $F(T)$. Furthermore, $F$ maps $\bsa$ onto the set
$Y:=\bigsetdef{S\in X}{\|S\|<1}\subset X$.

Denoting by $GX$ the invertible elements in $X$ and by $\cF X$ the Fredholm elements in $X$
(and similarly for Y) we find that
\begin{equation}
 F:\bigl((\cfsa,G\csa),d_R\bigr)\longrightarrow \bigl((\cF X,GX),\|\cdot\|\bigr)
\end{equation}
is a homeomorphism.

$F|\cB^\sa$ is a homeomorphism onto $Y$, too. Namely, by \cite[Prop. 1.5]{BooLesPhi:UFO}
the inverse of $F$ is given by $F^{-1}(S)=(1-S^2)^{-1/2}S$ and this is certainly
norm continuous on $Y$. Hence
\begin{equation}
 F:\bigl((\bfsa,G\bsa),\|\cdot\|\bigr)\longrightarrow \bigl((\cF Y,GY),\|\cdot\|\bigr)
\end{equation}
is a homeomorphism, too.

In sum, it suffices to prove that the inclusion
\begin{equation}
  \beta:=F\circ j \circ F^{-1}:(\cF Y,GY)\hookrightarrow (\cF X,GX)
\end{equation}
is a homotopy equivalence. Recall
that we are now dealing with sets of \emph{bounded} self--adjoint operators
which are equipped with the usual norm topology. Therefore, the map
\begin{equation}
    H:X\times [0,1/2]\longrightarrow X, (S,t)\mapsto S(I+S^2)^{-t}
\end{equation}
is trivially continuous. Moreover, it has the mapping properties
\begin{equation}\begin{split}
    H(GX\times [0,1/2])&\subset GX,\\
    H(Y\times [0,1/2])&\subset Y,\\
    H(X\times (0,1/2])&\subset Y,\\
    H(GY\times [0,1/2])&\subset GY.
\end{split}
\end{equation}
Put $g:(\cF X,GX)\to (\cF Y,GY), g:=H(\cdot,1/2)$. Then $g$ is a homotopy inverse
of $\beta$ since 
$H$ is a homotopy between $\id_{(\cF X,GX)}$ and $\beta\circ g$ and
the restriction of $H$ to $Y\times [0,1/2]$ is a homotopy between
$\id_{(\cF Y,GY)}$ and $g\circ \beta$.
\end{proof}
\begin{remark}\indent\par

\enum{1} Theorem \plref{uniqueness-Riesz} was observed in \cite[Sec. 3]{Nic:SFO}
without proof.

\enum{2}
We leave it to the reader to calculate the homotopy
inverse $F^{-1}\circ g\circ F$ of $j$ and the corresponding homotopy.
It is a tedious formula. Intuitively one would try the map $F$ itself
to be a homotopy inverse of $j$ and the homotopy to be
(same formula as $H$) $(T,s)\mapsto T(I+T^2)^{-s}, 0\le s\le 1/2$. 
However, for unbounded $T$ the operator $T(I+T^2)^{-s}$
is bounded for $s=1/2$ and unbounded for $s<1/2$ and 
proving continuity at $s=1/2$ seems to be tedious, though we did
not try very hard.
\end{remark}

As a consequence of Theorem \plref{uniqueness-Riesz} we note:

\begin{cor} Let $H$ be an infinite--dimensional separable complex Hilbert space.
Then $(\cfsa_*,d_R)$ is a classifying space for the $K^1$--functor.
Its homotopy groups are given by \eqref{G-I.2} and the uniqueness for the
spectral flow holds as in Theorem \plref{uniqueness-bounded}.
\end{cor}

Recall that in subsection \plref{secfivethree} we give a proof of Theorem
\plref{uniqueness-bounded} which is independent of \cite{AtiSin:ITS}. 

\subsection{Uniqueness of the spectral flow for the $d_W$--metric}

For completeness we state the uniqueness for the spectral flow in $(\cF\bsa,d_W)$.

\begin{theorem}
\label{uniqueness-dW} Let $H$ be an infinite--dimensional separable Hilbert
space and let $D$ be a fixed self--adjoint operator in $H$, $W:=\cD(D)$. Let
\[\mu:\Omega(\cF\bsa(W,H),G\bsa)\longrightarrow \Z\]
be a map which satisfies
\emph{Concatenation, Homotopy} and \emph{Normalization} in the following sense:
\begin{enumerate}
  \item[\relax] There is a rank one orthogonal projection $P$ with $\im P\subset W$
such that for all $A\in \bsa(W, H)$ with $(I-P)A(I-P)$ invertible we have
     \[\mu\bigl((tP+(I-P)A(I-P))_{-1/2\le t\le 1/2}\bigr)=1.\]
   \end{enumerate}
\end{theorem}

The \emph{Normalization} condition is slightly more complicated here. Superficially,
this is because we have formulated the theorem for paths in $\cF\bsa(W,H)$ instead
of $\cF_*\bsa(W,H)$. But this is not really the point. The problem is that we
do not even know whether $\cF_*\bsa(W,H)$ is path connected or not. If it is
then \emph{Normalization} can be formulated as in Theorem \plref{uniqueness-graph}.

Theorem \plref{uniqueness-dW} is basically 
due to Robbin and Salamon \cite{RobSal:SFM}, who
assumed additionally that $D$ has compact resolvent. The formulation in loc. cit.
is slightly different since they impose a \emph{Direct Sum} axiom.

The proof of Theorem \plref{uniqueness-dW} just follows along the lines of the
proof of Theorem \plref{uniqueness-bounded}: Assertions 1 and 2 just carry over
word by word. At first glance the family of unitaries chosen in the proof of Assertion 2
might be problematic. However, since $E(t)$ is finite--rank with image in $W$ one
sees that $U(t)$ can be chosen as a finite--rank perturbation of $I$ and such
that $U$ maps $W$ into itself. The proof of Assertion 3 uses that $(G\csa,d_G)$ is
path connected. Here this is taken care of by the stronger \emph{Normalization}
condition which allows to skip Assertion 3 and go directly to the ``Finish
of Proof'' on page \pageref{finishproof}. We leave the details to the reader.

\begin{appendix}
\section{Some estimates}\label{appendix}

In this appendix we collect a couple of operator estimates which
are basically well--known but for which references are hard to find.

%%%%%%%%%%%%%%%%%%%%%%%%%%%%%%%%%%%%%%%%%%%%%%%%%%%%%%%%%%%%%%%%%%%%%%%%%%%
% A functional analytic lemma

\begin{prop}\label{S-A.1} Let $H$ be a separable Hilbert space and let 
$T$ be an (unbounded) self--adjoint operator in $H$ with
bounded inverse. Furthermore, let $B$ be a symmetric operator in $H$
with $\cD(B)\supset\cD(T)$.

\begin{enumerate}
\item $T^{-1}BT$ is densely defined and $(T^ {-1}BT)^*=TBT^{-1}$.
\item If $T^{-1}BT$ or $TBT^{-1}$ is densely defined and bounded
then $B, TBT^{-1}$ and $T^{-1}BT$
are densely defined and bounded and we have $\|TBT^{-1}\|=\|T^{-1}BT\|$ and
$\|B\|\le \|T^{-1}BT\|$. 
\item If $T^{-1}B$ is bounded then so is $T^{-1/2}BT^{-1/2}$ and
$\|T^{-1/2}BT^{-1/2}\|\le \|T^{-1}B\|$.
\end{enumerate}
\end{prop}
Note that, by definition,
\begin{equation}\label{G-A0}\begin{split}
 \cD(TBT^{-1})&=\bigsetdef{x\in H}{BT^{-1}x\in\cD(T)}\\
 \cD(T^{-1}BT)&=\bigsetdef{x\in\cD(T)}{Tx\in\cD(B)}.
			    \end{split}
\end{equation}
In \enum{1} it is \emph{not} claimed that $\cD(TBT^{-1})$
is dense. Hence if $T^{-1}BT$ is bounded then it \emph{follows} that
$TBT^{-1}$ is defined on $H$ and bounded. Note that by \enum{1} the operator
$TBT^{-1}$ is always closed.

\begin{proof} This is basically a consequence of complex interpolation
theory (Taylor \cite[Sec. 4.2]{Tay:PDEI}) but we prefer to give a direct elementary
proof here. 

We first note that \enum{3} follows from \enum{2}: namely
the operator $X:=T^{-1/2}BT^{-1/2}$ is symmetric on $\cD(T^{1/2})$ and
$T^{-1/2}XT^{1/2}$ is densely defined and bounded. 
Now apply \enum{2} with $B=X$ and $T^{1/2}$
instead of $T$.

To prove \enum{1} we note that certainly $\cD(T^{-1}BT)=\bigsetdef{x\in
  \cD(T)}{Tx\in\cD(B)}\supset\cD(T^2)$. Since
$T$ is self--adjoint $\cD(T^2)$ is a core for $T$ (i.e. $\cD(T^2)$ is dense
in $\cD(T)$ with respect to the graph norm), in particular it is dense in $H$. 
Thus $T^{-1}BT$ is densely defined. 

Now consider $x\in\cD(TBT^{-1})$ and $y\in\cD(T^{-1}BT)$. Then, by \eqref{G-A0},
$y\in\cD(T), Ty\in\cD(B)$
and $T^{-1}x\in\cD(B), BT^{-1}\in\cD(T)$ and consequently,
\begin{equation}
  \scalar{TBT^{-1}x}{y}=\scalar{BT^{-1}x}{Ty}=\scalar{T^{-1}x}{BTy}=\scalar{x}{T^{-1}BTy}.
\end{equation}
This shows $TBT^{-1}\subset(T^{-1}BT)^*$.

To show the converse inclusion, let us consider $x\in\cD((T^{-1}BT)^*)$ and
$y\in\cD(T^2)\subset\cD(T^{-1}BT)$. Carefully checking domains we find
\begin{equation}\label{G-A5}
  \scalar{(T^{-1}BT)^*x}{y}=\scalar{x}{T^{-1}BTy}=\scalar{T^{-1}x}{BTy}=\scalar{BT^{-1}x}{Ty}. 
\end{equation}
As noted above, $\cD(T^2)$ is dense in $\cD(T)$ with respect to the graph norm.
Hence the equality
\begin{equation}
    \scalar{(T^{-1}BT)^*x}{y}=\scalar{BT^{-1}x}{Ty}
\end{equation}
holds for all $y\in\cD(T)$. But this means that $BT^{-1}x\in\cD(T)$ and
$TBT^{-1}x=(T^{-1}BT)^*x$ proving $TBT^{-1}\supset(T^{-1}BT)^*$.

To prove \enum{2} assume that the densely defined operator 
$T^{-1}BT$ is bounded. Then its adjoint $TBT^{-1}=(T^{-1}BT)^*$ is densely defined and
bounded, too.

If $TBT^{-1}$ is densely defined and bounded on $H$ 
then we infer from $TBT^{-1}=(T^{-1}BT)^*$ that
$T^{-1}BT$ is closable and $\overline{T^{-1}BT}=(TBT^{-1})^*$. Hence $T^{-1}BT$
is (densely defined and) bounded.

It is now clear that if $TBT^{-1}$ and $T^{-1}BT$ are bounded that then
$\|TBT^{-1}\|=\|T^{-1}BT\|$.

Now suppose that they are both bounded and pick $x,y\in\cD(T^2)$ 
and consider the analytic function
\begin{equation}
     f(z):=\scalar{T^{2z-1}BT^{1-2z}x}{y},\quad 0< \Re z< 1.
\end{equation}
Since $x,y\in\cD(T^2)$ it is straightforward
to check that $f$ is bounded and continuous
on the vertical strip $\bigsetdef{z\in\C}{0\le \Re z\le 1}$. Moreover,
we have for $z=it, t\in\R, $
\begin{equation}
   |f(z)|=|\scalar{T^{-1}BTT^{-2it}x}{T^{-2it}y}|\le \|T^{-1}BT\|\,\|x\|\,
   \|y\|,
\end{equation}
and similarly for $z=1+it,t\in\R,$
\begin{equation}
   |f(z)|=|\scalar{T^{-2it}x}{T^{-1}BT^{1-2it}y}|\le 
         \|T^{-1}BT\|\,\|x\|\,\|y\|.
\end{equation}
Hence by Hadamard's three line theorem (Rudin \cite[Thm. 12.8]{Rud:RCA}) we find
$|f(z)|\le \|T^{-1}BT\|\,\|x\|\,\|y\|$ for $0\le \Re z\le 1$. In particular
we have for $z=1/2$
\begin{equation}  |\scalar{Bx}{y}|\le \|T^{-1}BT\|\,\|x\|\,\|y\|,\quad x,y\in\cD(T^2).
\end{equation}
Since $\cD(T^2)$ is dense in $H$ we reach the conclusion.
\end{proof}
%%%%%%%%%%%%%%%%%%%%%%%%%%%%%%%%%%%%%%%%%%%%%%%%%%%%%%%%%%%%%%%%%%%%%%%%%

\begin{prop}\label{S-A.2} Let $T\in\csa, S\in\bsa$ with
$(T+i)S(T+i)^{-1}$ densely defined and bounded. Then, for $0\le \ga\le 2, -1\le\beta\le 1,$
$\ga+\beta<2$ we have the norm estimate
\[\begin{split}
  \Bigl\| |T+S+i|^\ga&\bigl(|T+S+i|^{-1}-|T+i|^{-1}\bigr)|T+i|^\beta\Bigr\|\\
&\le C(\ga,\beta) \Bigl(\bigl\|(T+i)S(T+i)^{-1}\bigr\|+\|S^2(T+i)^{-1}\|\Bigr).
  \end{split}
\]
\end{prop}
\begin{remark}\label{S-A.3} For $(T+S+i)^{-1}-(T+i)^{-1}$ in place of 
$|T+S+i|^{-1}-|T+i|^{-1}$ the estimate follows easily from the
resolvent equation
\begin{equation}
    (T+S+i)^{-1}-(T+i)^{-1}=-(T+S+i)^{-1}S(T+i)^{-1}
\end{equation}
and complex interpolation theory.

To deal with the operator absolute value recall that 
for any non--negative invertible operator $A$ in $H$ one has
\begin{equation}\label{G-A1}
    A^{-1/2}=\frac{2}{\pi}\int_0^\infty (A+x^2)^{-1}dx.
\end{equation}
\end{remark}
\begin{proof}
Note that $(T+i)S(T+i)^{-1}$ bounded means that $S$ maps
$\cD(T)$ continuously into $\cD(T)$. Hence $\cD((T+S)^2)\supset\cD(T^2)$
and thus $(T+S)^2-T^2=(T+S)S+ST$ on $\cD(T^2)$. 
Thus we have for $x\ge 0$ the resolvent identity
\begin{equation}\begin{split}
   (I+&(T+S)^2+x^2)^{-1}-(I+T^2+x^2)^{-1}\\
     &= -(I+(T+S)^2+x^2)^{-1}\bigl((T+S)S+ST\bigr)(I+T^2+x^2)^{-1}\\
     &=:\cI(T,S,x),
                \end{split}\label{G-A2}
\end{equation}
and in view of \eqref{G-A1} we find
\begin{equation}\begin{split}
    \bigl( &|T+S+i|^{-1}-|T+i|^{-1}\bigr)\\
     &=-\frac{2}{\pi}\int_0^\infty \cI(T,S,x) dx.
                \end{split}\label{G-A.8}
\end{equation}
Next we estimate the integrand of \eqref{G-A.8}
\begin{equation}\begin{split}
    \Bigl\| |T+&S+i|^\ga\cI(T,S,x)|T+i|^\beta\Bigr\|\\
     &\le  \Bigl\| |T+S+i|^\ga(I+(T+S)^2+x^2)^{-1}\Bigr\|\cdot\ldots\\
        &\qquad\cdot
          \bigl\|(T+S)S|T+i|^{-1}\bigr\| 
           \Bigl\||T+i|(I+T^2+x^2)^{-1}|T+i|^\beta\Bigr\|\\
     &\quad +\Bigl\| |T+S+i|^\ga(I+(T+S)^2+x^2)^{-1}\Bigr\|\cdot\ldots\\
        &\qquad \cdot \|S\| \Bigl\|T(I+T^2+x^2)^{-1}|T+i|^\beta\Bigr\|\\
     &\le C(\ga,\beta)(1+x^2)^{(\ga+\beta-3)/2}
       \bigl(\|S\|+\|(T+S)S|T+i|^{-1}\|\bigr)\\
     &\le C(\ga,\beta)(1+x^2)^{(\ga+\beta-3)/2}
       \bigl(\||T+i|S|T+i|^{-1}\|+\|(S^2|T+i|^{-1}\|\bigr).
                \end{split}\label{G-A4}
\end{equation}
In the last inequality we have used Proposition \plref{S-A.1}.

If $\ga+\beta<2$ we may integrate \eqref{G-A4} and reach the conclusion.
\end{proof}
\end{appendix}

% depending on style, choose one of the following
%\newcommand{\Toappear}{To appear in}  
\newcommand{\Toappear}{to appear in}

% if you want to use bibtex, uncomment the following two lines
% and comment out the next input command
%\bibliography{mlabbr,papers2003,books2003}
%\bibliographystyle{amsalpha}
\providecommand{\bysame}{\leavevmode\hbox to3em{\hrulefill}\thinspace}
\providecommand{\MR}{\relax\ifhmode\unskip\space\fi MR }
\providecommand{\MRhref}[2]{%
  \href{http://www.ams.org/mathscinet-getitem?mr=#1}{#2}
}
\providecommand{\href}[2]{#2}

\newcommand{\foottext}[1]{{%
  \renewcommand{\thefootnote}{}\footnotetext{#1}}}

%\foottext{Received by the editors January 29, 2004; revised April 29, 2004.}

\end{document}